%% file: 00000_main.tex
\documentclass[final,1p,times]{elsarticle}
\input{00_packages}

\input{00_dougcommands}

\journal{Computer \& Fluids}

\title{Full-Space Approach to Aerodynamic Shape Optimization}
\begin{document}

\begin{frontmatter}

\author{Doug Shi-Dong\corref{cor1}\fnref{fn1}}
\ead{doug.shi-dong@mail.mcgill.ca}
\cortext[cor1]{Corresponding author}
\fntext[fn1]{Doctoral Candidate}

\author{Siva Nadarajah\fnref{fn2}}
\ead{siva.nadarajah@mcgill.ca}
\fntext[fn2]{Associate Professor}

\address{Department of Mechanical Engineering, McGill University, Montreal, Quebec, Canada}

\input{1_abstract}


\end{frontmatter}

\input{2_introduction}
\input{3_optimization_problem}
\input{4_flow_constraints}

\input{5_design_parametrization}
\input{6_meshmover}
\input{7_cost_breakdown}
\input{7_results}
\input{9_conclusionandacknowledgements}


%\section*{References}

\clearpage
\bibliographystyle{elsarticle-num} 
\bibliography{000_references}

\input{9z_appendix}
\typeout{get arXiv to do 4 passes: Label(s) may have changed. Rerun}

\end{document}

%% file: 00_packages.tex
\usepackage{mathtools}

\usepackage[svgnames,table]{xcolor}% note the table option
\usepackage{booktabs}

\usepackage{amssymb}
\usepackage{amsmath}
\usepackage{subfig}
\usepackage{hyperref}
\usepackage{wrapfig}
\usepackage{graphicx}
\usepackage[makeroom]{cancel} % cross out a term
\usepackage{soul} %striketrough
\usepackage{siunitx} % Table formatting
\usepackage{fp} % Do some simple arithmetics

\usepackage{booktabs}
\usepackage{pifont}% http://ctan.org/pkg/pifont Checkmarks
\usepackage{cleveref} % Allows \cref{eq:eq1, eq:eq2, eq:eq3} to look good

\usepackage{pgfplots}
\pgfplotsset{compat=1.8}
\usepgfplotslibrary{colorbrewer,patchplots}

\usepackage{tikz} % Flowchart
\usetikzlibrary{shapes.geometric, fit, arrows, decorations.markings, shapes, calc, positioning}
\usepackage{color}
\usepackage[normalem]{ulem}

%% file: 00_dougcommands.tex
\newcommand{\Eminus}[1]{\text{\sc{e}-}#1}
\newcommand{\oneEm}[1]{1.0\Eminus{#1}}

\newcommand{\mbf}[1]{\mathbf{#1}}
\newcommand{\mbs}[1]{\boldsymbol{#1}}

\renewcommand{\d}{\mathop{}\!\mathrm{d}}
\newcommand{\del}[2]{
\dfrac{\partial#1}{\partial#2}}

\newcommand{\Dt}[2]{\dfrac{\mathrm{d}#1}{\mathrm{d}#2}}
\newcommand{\DDt}[2]{\dfrac{\mathrm{d}^2#1}{\mathrm{d}#2^2}}

\newcommand{\bx}[0]{\mbf{x}}

\newcommand{\bzero}[0]{\mbf{0}}

\newcommand{\bOmega}[0]{\mbs{\Omega}}

\newcommand{\identity}[0]{\mbf{I}}

\newcommand{\real}[0]{\mathbb{R}}
\newcommand{\domain}[0]{\bOmega}
\newcommand{\boundary}[0]{{\partial \bOmega}}

\newcommand{\objective}[0]{\mathcal{F}}

\newcommand{\des}[0]{{\boldsymbol{\alpha}}}

\newcommand{\KKTmat}[0]{\mbs{\mathcal{K}}}
\newcommand{\PFour}[0]{\mbf{P}_4}
\newcommand{\PFourA}[0]{\Tilde{\mbf{P}}_4}
\newcommand{\PTwo}[0]{\mbf{P}_2}
\newcommand{\PTwoA}[0]{\Tilde{\mbf{P}}_2}

\newcommand{\augpenalty}[0]{\mu}
\newcommand{\lagrangian}[0]{\mathcal{L}}
\newcommand{\dual}[0]{\mbs{\lambda}}

\newcommand{\approxHessian}[0]{\mbf{B}_{z}}
\newcommand{\search}[1]{\mbf{p}_{#1}}

\newcommand{\residual}{\mathcal{R}}
\newcommand{\residuals}{\mbs{\residual}}
\newcommand{\state}{\mbf{u}}
\newcommand{\impstate}{\mbf{w}}

\newcommand{\xvol}[0]{\mbf{x_{vol}}}
\newcommand{\xsurf}[0]{\mbf{x_{surf}}}

\newcommand{\ndes}[0]{n}
\newcommand{\nstate}[0]{m}
\usepackage{stmaryrd}

\usepackage{mleftright}
\mleftright % make \left & \right behave like \mleft & \mright 

\definecolor{bostonuniversityred}{rgb}{0.8, 0.0, 0.0}

\usepackage{color, soul}
\sethlcolor{yellow}
\makeatletter
\def\SOUL@hlpreamble{%
    \setul{\dp\strutbox}{\dimexpr\ht\strutbox+\dp\strutbox\relax}%
    \let\SOUL@stcolor\SOUL@hlcolor
    \SOUL@stpreamble
}
\makeatother
\DeclareRobustCommand{\corr}[1]{%
  \ifmmode\text{\hl{$#1$}}\else\hl{#1}\fi
}
\usepackage{array,tabularx}

%% file: 1_abstract.tex
\begin{abstract}
Aerodynamic shape optimization (ASO) involves finding an optimal surface while constraining a set of nonlinear partial differential equations (PDE).
The conventional approaches use quasi-Newton methods operating in the reduced-space, where the PDE constraints are eliminated at each design step by decoupling the flow solver from the optimizer.
Conversely, the full-space Lagrange-Newton-Krylov-Schur (LNKS) approach couples the design and flow iteration by simultaneously minimizing the objective function and improving feasibility of the PDE constraints, which requires less iterations of the forward problem.
Additionally, the use of second-order information leads to a number of design iterations independent of the number of control variables.
We discuss the necessary ingredients to build an efficient LNKS ASO framework as well as the intricacies of their implementation.
The LNKS approach is then compared to reduced-space approaches on a benchmark two-dimensional test case using a high-order discontinuous Galerkin method to discretize the PDE constraint.
\end{abstract}

%% file: 2_introduction.tex
\section{Introduction}
%The field of Aerodynamic shape optimization (ASO) has seen an increase in both
%the level of complexity of the geometry leading to even increasing computational grids
%and the number of control variables to design the geometry {\color{red} [need citations here]}.
The field of aerodynamic shape optimization (ASO) is at a stage where complete aircraft~\cite{Reuther1996,Palacios2015,Gagnon2016,Chen2016} and engine~\cite{Luo2010, Walther2015} configurations can be designed through numerical optimization.
As a result, increasingly complex geometries lead to ever increasing computational grids.
%has seen an increase in both
%the level of complexity of the geometry leading to even increasing computational grids
%and the number of control variables to design the geometry {\color{red} [need citations here]}.
While the number of control variables for shape optimization is currently in the order of tens or hundreds,
% and will likely not exceed in the order of a few thousands,
the need for higher fidelity simulations is requiring the convergence of larger than ever
discretized systems of partial differential equations (PDE).
It is therefore advantageous to require a valid flow solution as infrequently as possible throughout the design process.
%It is therefore in the optimizer's interest to require a valid flow solution as infrequently as possible.
%The solution of a linear system of the size of the flow constraints will be deemed as one unit of work.

There are two ways to reduce the number of linear systems to be solved: decrease the number of design cycles or lower the number of linear systems to be solved within each design cycle.
While current reduced-space quasi-Newton optimization frameworks can evaluate first-order sensitivities independently of the problem size through the adjoint method~\cite{Jameson1988}, the number of design cycles required to converge still increases proportionally to the number of control variables as more quasi-Newton updates are needed to correctly approximate the Hessian.
Additionally, the reduced-space approach typically performs multiple iterations of the flow solver until full convergence of the PDE constraint, although
tolerances on the flow residuals have been derived to maintain the same Q-convergence rate of the optimization problem~\cite{Brown2017}.

Using Newton iterations within the reduced-space allows the optimizer to converge independently of the control variables size at the detriment of a higher cost within each cycle.
While second-order adjoints~\cite{Sherman1996,Ghate2007,Papadimitriou2008} allow the Hessian to be computed at a cost proportional to the number of control variables, we have basically shifted the work dependence on the control dimensions from the outer iteration to the inner iteration.
Hessian approximations~\cite{ShiDong2018,Arian1999,Schmidt2010} have been formulated to reduce its prohibitive cost and have demonstrated significant convergence improvement.
Another possible solution is to solve the reduced system using Krylov subspace methods, which only require Hessian-vector products, at the cost of two linear solves per Krylov iteration~\cite{Heinkenschloss1999,Akcelik2002,Hicken2014}.
However, the number of Krylov iterations required would scale with the size of the reduced system.
%When an inexact matrix-free Krylov approach \cite{Heinkenschloss1999,Akcelik2002,Hicken2014} is used, each Hessian-vector product requires two linear solves.

%While current quasi-Newton optimization frameworks can evaluate first-order sensitivities independently of the problem size through the adjoint method \cite{Jameson1988}, the number of design iterations required to converge still increases proportionally.
%Using Newton steps is therefore attractive to attain convergence independently of the number of design variables.
%However, each Newton step is associated with a higher cost.
%In the reduced-space approach, the reduced-Hessian requires as many linear solves as the number of design variables.
%When an inexact matrix-free Krylov approach \cite{Heinkenschloss1999,Akcelik2002,Hicken2014} is used, each Hessian-vector product requires two linear solves.

The full-space Lagrange-Newton-Krylov-Schur (LNKS) from Biros and Ghattas~\cite{Biros2005a,Biros2005} framework hopes to achieve convergence independent of problem size if the user can provide optimal preconditioners of the forward problem, dual problem, and reduced-Hessian.
%Their proposed preconditioners
However, its implementation on complex problems have been challenging since the optimization framework must be able to store and manipulate parallelized discretizations rather than letting a blackbox solver handle the parallelization.
For example, popular packages in the ASO community, such as SNOPT~\cite{Gill2005}, NLPQLP~\cite{Schittkowski2009}, or IPOpt~\cite{Waechter2005} cannot be used to implement LNKS.
Apart from a quasi-one dimensional implementation~\cite{Hicken2013,Hicken2014}, the use of a Newton-based full-space approach has yet to be achieved in ASO on a larger scale.
%needs to parallelize its storage and perform parallel matrix-vector products large-scale discretizations.
%Additionally, the previously proposed full-space preconditioner build on our knowledge of preconditioning the reduced-space systems.
%Biros and Ghattas \cite{Biros2005a, Biros2005} use \texttt{PETSc} \cite{petsc-efficient} to handle the PDE discretization and optimization framework.

In the present work, we define the different existing frameworks, and the cost associated with their sensitivities requirements.
A description of the ingredients necessary to implement LNKS within an ASO framework is then provided.
The resulting preconditioned systems in~\ref{app:preconditioner_inverses} and re-derived in~\ref{app:precond_proof} amend the original work~\cite{Biros2005a}. Despite the differences, the eigenvalue spectra of the preconditioned systems shown in~\ref{app:precond_eig} coincide with the original preconditioner proposed in~\cite{Biros2005a}, and therefore do not affect the conclusions in~\cite{Biros2005a}.

We perform a dimensional analysis for an increasing number of control variables, and an increasing number of state variables to showcase the scalability of the various algorithms.
Finally, we compare the full-space approach against reduced BFGS and reduced full-space methods similarly to~\cite{Hicken2013} on an inverse design problem.

%% file: 3_optimization_problem.tex
\section{Optimization formulation}
\label{sec:optimization_formulation}
\subsection{Problem statement}
The ASO problem formulation can be stated as 
\begin{equation}
    \min_{\state,\des} \objective(\state,\des) 
    \quad
    \text{subject to}
    \quad
    \residuals(\state,\des) = \bzero
    ,
\label{eq:optimization_problem}
\end{equation}
where
$\des \in \real^\ndes$ are the control variables
, $\state \in \real^\nstate$ are flow state variables
, $\objective: \real^{\nstate} \times \real^{\ndes} \to \real$ is the objective function, and $\residuals: \real^{\nstate} \times \real^{\ndes} \to \real^m$ are the PDE constraints relating the geometry and flow solution.

The reduced and full-space approaches are implemented using the Trilinos' Rapid Optimization Library (\texttt{ROL}) \cite{rol-website} which provides a modular large-scale parallel optimization framework.

\input{3a_full_space}
\input{3b_reduced_space}

\input{3c_precond_full_space}

%% file: 3a_full_space.tex
\subsection{Full-space}
The full-space approach directly tackles the original problem in Eq.~\eqref{eq:optimization_problem}.
%\begin{equation}
%    \min_{\state,\des} \objective(\state,\des) 
%    \quad
%    \text{subject to}
%    \quad
%    \residuals(\state,\des) = \bzero
%    .
%\end{equation}
To recover the optimality conditions, we form the Lagrangian function $\lagrangian$ as
\begin{equation}
    \lagrangian(\state,\des) \coloneqq \objective(\state,\des) + \dual^T \residuals(\state,\des)
    ,
\end{equation}
where $\dual \in \real^\nstate$ are the adjoint variables.
The first-order necessary optimality conditions, also known as the Karush-Kuhn-Tucker (KKT) conditions, are given by
\begin{subequations}
%\begin{align}
%    \del{\lagrangian(\state,\des)}{\state} &= \del{\objective(\state,\des)}{\state} + \dual^T \del{\residuals(\state,\des)}{\state} = \bzero,
%    \\
%    \del{\lagrangian(\state,\des)}{\des} &= \del{\objective(\state,\des)}{\des} + \dual^T \del{\residuals(\state,\des)}{\des} = \bzero,
%    \\
%    \del{\lagrangian(\state,\des)}{\dual} &= \residuals(\state,\des) = \bzero.
%\end{align}
\begin{align}
    {\lagrangian_\state(\state,\des)} &= {\objective_\state(\state,\des)} + \dual^T {\residuals_\state(\state,\des)} = \bzero,
    \label{eq:kkt_conditions_state}
    \\
    {\lagrangian_\des(\state,\des)} &= {\objective_\des(\state,\des)} + \dual^T {\residuals_\des(\state,\des)} = \bzero,
    \label{eq:kkt_conditions_control}
    \\
    {\lagrangian_{\dual}(\state,\des)} &= \residuals(\state,\des) = \bzero.
    \label{eq:kkt_conditions_adjoint}
\end{align}
\label{eq:kkt_conditions}
\end{subequations}
The Newton iteration represented by the KKT system 
%solving those equations requires the solution to the KKT system
\begin{equation}
\KKTmat \search{} =
    \begin{bmatrix}
    \lagrangian_{\state \state} & \lagrangian_{\state \des} &  \residuals_{\state}^T
    \\
    \lagrangian_{\des \state} & \lagrangian_{\des \des} &  \residuals_{\des}^T
    \\
    \residuals_{\state} & \residuals_{\des} & \bzero
    \end{bmatrix}
    \begin{bmatrix}
    \search{\state} \\ \search{\des} \\  \search{\dual}
    \end{bmatrix}
    =
    -
    \begin{bmatrix}
    \objective_\state + \dual^T \residuals_\state
    \\
    \objective_\des + \dual^T \residuals_\des
    \\
    \residuals
    \end{bmatrix}
    ,
\label{eq:kkt_system}
\end{equation}
%\begin{equation}
%    \mbf{K} \mbf{p} = \mbf{g}
%    ,
%\label{eq:abbr_kkt_system}
%\end{equation}
where $\KKTmat$ represents the KKT matrix and $\search{}$ is composed of the search directions $\search{\state} = \state^{i+1} - \state^{i}$, $\search{\des} = \des^{i+1} - \des^{i}$, and $\search{\dual} = \dual^{i+1} - \dual^{i}$.
Its assembly requires second-order sensitivities of the objective function and residuals to form the Hessian of the Lagrangian.

The dimensions of the KKT matrix ($2\nstate+\ndes$) on the left-hand side of Eq.~\eqref{eq:kkt_system} precludes direct factorization methods and is therefore solved using Krylov subspace methods such as \mbox{GMRES}~\cite{Saad1986}.
Since it is a Newton-based method, the number of times this system needs to be solved is independent of the problem size.
Furthermore, its solution only necessitates matrix-vector products, which involves no inversion of the forward or backward problem.

\subsubsection{Full-space implementation}
The search direction $\search{}$ is obtained by solving the KKT system Eq.~\eqref{eq:kkt_system} with FGMRES.
A relative tolerance of \oneEm{6} is used to determine convergence of the linear system.
The various preconditioners proposed by~\cite{Biros2005a} are discussed in Sec.~\ref{sec:full_space_precond}.
A backtracking line-search also complements the algorithm for globalization purposes, where the augmented Lagrangian is used as the merit function as discussed in~\cite{Biros2005}
\begin{equation}
    \lagrangian^{+} \coloneqq \objective + \dual^T \residuals + \frac{\augpenalty}{2} \residuals^T \residuals
    = \lagrangian + \frac{\augpenalty}{2} \residuals^T \residuals,
\end{equation}
with the penalty parameter $\augpenalty$.

Our full-space implementation starts with a converged flow solution such that the flow residual is zero.
While the full-space method does not postulate a feasible initial solution, it is performing Newton steps from the perspective of the steady-state flow solver.
Therefore, starting from a free-stream solution resulted in inappropriate search directions, while a converged flow starts in the Newton ball of convergence.

%The penalty parameter defined by Biros and Ghattas $\augpenalty$ is inversely proportional to the residual.
However, the penalty parameter defined by Biros and Ghattas~\cite{Biros2005}, 
\begin{equation}
    \augpenalty_{\text{BG}} = 
    \frac{
        \search{\state}^T \lagrangian_{\state} 
        + \search{\des}^T \lagrangian_{\des} 
        + \search{\dual}^T \lagrangian_{\dual} + \delta
    }
    {\residuals^T \left(\residuals_{\state}\search{\state} + \residuals_{\des}\search{\des} \right)
    }
\end{equation} is inversely proportional to the residual, $\residuals$. Therefore, the initial converged flow solution, heavily penalized constraint violations, resulting in small step-lengths. Instead, we defined the penalty parameter to be inversely proportional to the norm of the Lagrangian gradient with respect to the control variables
\begin{equation}
    \augpenalty = \frac{c_1}{\|\lagrangian_{\des}\|_2},
\end{equation}
such that constraints are enforced to a greater extent as the optimization converges.
A factor $c_1 = 0.01$ is used to scale the penalty parameter.

Previous one-shot methods leverage PDE solver technologies by using well-established explicit time-stepping schemes into their algorithm~\cite{Gauger2008_OptimalDesignOneShot}.
Since the current full-space approach is solved implicitly, pseudo-transient continuation~\cite{Kelley1998} is used to increase the robustness of the algorithm.
For steady-state flows, it is typical to cast the problem of solving
\begin{equation}
    \residuals(\state,\des) = \bzero
\end{equation}
into
\begin{equation}
    \mbs{\mathcal{M}}\del{\state}{\tau} + \residuals(\state,\des) = \bzero
,
\end{equation}
where $\mbs{\mathcal{M}}$ represents the mass matrix and $\tau$ the pseudo-time.
A backward Euler temporal discretization leads to 
\begin{equation}
    \left(\frac{\mbs{\mathcal{M}}}{\Delta \tau} + \residuals_{\state} \right) \left(\state^{i+1} - \state^{i} \right) = - \residuals
,
\label{eq:backward_euler}
\end{equation}
which is analogous to solving the third row of Eq.~\eqref{eq:kkt_system} if $\Delta \tau \to \infty$ and $\search{\des} \to \bzero$
\begin{equation}
    \residuals_{\state} \search{\state} = - \residuals - \residuals_{\des} \search{\des}
.
\end{equation}
%In fact, a large pseudo-time step $\Delta \tau$ and zero search direction $\search{\des}$ gives the
%An additional mass-matrix from the discretization is scaled by a a pseudo-timestep such that the flow 
The Jacobian $\residuals_{\state}$ is replaced by
\begin{equation}
    \residuals_{\state}^{+} = \residuals_{\state} + \frac{\mbs{\mathcal{M}}}{\Delta \tau} 
,
\end{equation}
where the pseudo-time step inversely scales with the square of the Lagrangian gradient norm
\begin{equation}
    \Delta \tau = \frac{1}{\|\residuals \|^2_2}.
\end{equation}
While it is initially equivalent to solving a different set of constraints
\begin{equation}
    \residuals^{+} =  \frac{\mbs{\mathcal{M}}}{\Delta \tau} \state + \residuals
    ,
\end{equation}
the quadratically increasing time-step ensures that the continuation term disappears.

%A value of $c_2 = 1.0/\|\lagrangian_{\des}\|$ was chosen for this work.

%Just like in~\cite{Gauger2008_OptimalDesignOneShot}, replacing the flow Jacobian allows the user to leveragemore sophisticated 
%where $\mbs{\mathcal{M}}$ represents the mass matrix analogous to solving the backward Euler

%% file: 3b_reduced_space.tex
\subsection{Reduced-space}
Reduced-space approaches eliminate the constraints from the optimization problem by letting the flow solver satisfy the constraints at all times.
%from Eq. \eqref{eq:kkt_conditions_adjoint} at all times.
Assuming that the flow Jacobian $\residuals_{\state} \in \real^{\nstate \times \nstate}$ is invertible, the implicit function theorem states that there exists a function $\impstate : \real^\ndes \to \real^\nstate$ such that
\begin{equation}
    \impstate(\des) \equiv \state
    \quad \text{and} \quad
    \residuals(\des, \impstate(\des)) = \bzero.
\end{equation}
The residual's total derivatives are then given by
\begin{equation}
    \Dt{\residuals}{\des} = \residuals_\state \impstate_\des + \residuals_\des = \bzero
    \label{eq:residual_total_first_derivative}
    ,
\end{equation}
where the sensitivity
\begin{equation}
    \impstate_\des = - \residuals_\state^{-1} \residuals_\des
    \label{eq:first_order_flow_sensitivity}
    ,
\end{equation}
requires $\ndes$ linear system solves.

In the reduced approach, a flow solver becomes responsible to ensure that the residual constraints are satisfied at every design cycle.
For a nonlinear PDE, this involves multiple fixed-point iterations of the forward problem, which involves solving the flow Jacobian linear system multiple times until convergence.

The optimization problem can then be recast as the reduced formulation
\begin{equation}
    \min_{\des} \objective(\des, \impstate(\des)).
\end{equation}

\subsubsection{Reduced gradient}
The reduced gradient $\mbf{g}_z$ of the objective function with respect to the control variables is given by
\begin{equation}
    \mbf{g}_z 
    = \Dt{\objective}{\des} 
    = \objective_\des + \objective_\state \impstate_\des
    = \objective_\des - \objective_\state \residuals_\state^{-1} \residuals_\des
    = \objective_\des + \dual^T \residuals_\des
    ,
\label{eq:reduced_gradient}
\end{equation}
where the adjoint equation
\begin{equation}
    \residuals_\state^{T} \dual = -\objective_\state^T
    ,
\label{eq:adjoint_equation}
\end{equation}
requires a single linear system solve to yield the dual state $\dual \in \real^\nstate$ .
%The reduced space therefore satisfies the KKT condition Eq

%The dual variable satisfying Eq. \eqref{eq:adjoint_equation} leads to a zero adjoint residual
%\begin{equation}
%    \residuals^{\dual} \coloneqq 
%    \objective_\state +  \dual^{T} \residuals_\state = \bzero,
%\label{eq:adjoint_residual}
%\end{equation}
%which will be used to derive the reduced Hessian matrix-vector product.

\subsubsection{Reduced Hessian}
Using a similar approach as Eqs. \eqref{eq:residual_total_first_derivative} and \eqref{eq:first_order_flow_sensitivity}, we can derive the second-order residual and flow sensitivities
\begin{equation}
    \DDt{\residuals}{\des} 
    = \impstate_\des^T \residuals_{\state \state} \impstate_\des
    + \residuals_{\des \state} \impstate_\des
    + \impstate_\des^T \residuals_{\state \des}
    + \residuals_{\state} \impstate_{\des \des}
    + \residuals_{\des \des}
    = \bzero
    ,
\end{equation}
such that
\begin{equation}
    \impstate_{\des \des} = 
    - \residuals_{\state}^{-1}
    \left( 
        \impstate_\des^T \residuals_{\state \state} \impstate_\des
        + \residuals_{\des \state} \impstate_\des
        + \impstate_\des^T \residuals_{\state \des}
        + \residuals_{\des \des}
    \right)
    .
\end{equation}
The second-order total derivatives of the objective function is then given by
\begin{equation}
\begin{split}
    \DDt{\objective}{\des} 
    &= \objective_{\des \des}
    + \objective_{\des \state} \impstate_\des
    + \impstate_\des^T \objective_{\state \des} 
    + \impstate_\des^T \objective_{\state \state} \impstate_\des
    + \objective_\state \impstate_{\des \des}
% Replace wxx
    \\&=
    \objective_{\des \des}
    + \objective_{\des \state} \impstate_\des
    + \impstate_\des^T \objective_{\state \des} 
    + \impstate_\des^T \objective_{\state \state} \impstate_\des
    % preprint newline
    \\&\quad
    - \objective_\state \residuals_{\state}^{-1}
    \left( 
        \impstate_\des^T \residuals_{\state \state} \impstate_\des
        + \residuals_{\des \state} \impstate_\des
        + \impstate_\des^T \residuals_{\state \des}
        + \residuals_{\des \des}
    \right)
% Use adjoint
    \\&=
    \objective_{\des \des}
    + \objective_{\des \state} \impstate_\des
    + \impstate_\des^T \objective_{\state \des} 
    + \impstate_\des^T \objective_{\state \state} \impstate_\des
    % preprint newline
    \\&\quad
    + \dual^T
    \left( 
        \impstate_\des^T \residuals_{\state \state} \impstate_\des
        + \residuals_{\des \state} \impstate_\des
        + \impstate_\des^T \residuals_{\state \des}
        + \residuals_{\des \des}
    \right)
% Regroup terms
    \\&=
    \left( \objective_{\des \des} + \dual^T \residuals_{\des \des} \right)
    + \left ( 
        \objective_{\des \state} + \dual^T \residuals_{\des \state} 
      \right) \impstate_\des
    % preprint newline
    \\&\quad
    + \impstate_\des^T
      \left( \objective_{\state \des}  + \dual^T \residuals_{\state \des}
      \right)
    + \impstate_\des^T
      \left( \objective_{\state \state} + \dual^T \residuals_{\state \state} \right)
      \impstate_\des
    \\&=
    \lagrangian_{\des \des}
    + \lagrangian_{\des \state} \impstate_\des
    + \impstate_\des^T \lagrangian_{\state \des}
    + \impstate_\des^T \lagrangian_{\state \state}
      \impstate_\des,
\end{split}
\end{equation}
which requires $\ndes$ additional linear solves to obtain $\impstate_\des$, assuming that $\dual$ has been pre-computed for the gradient.
For ASO, the resulting matrix is usually small enough to directly be factorized.
However, the  $\ndes$ additional linear solves quickly become prohibitive and some approximations can reduce its cost \cite{ShiDong2018}.

%The same derivation can be obtained by choosing the appropriate constraint null-space and range-space basis as shown in \cite{Biros2005a,Nocedal2006}.
The reduced Hessian matrix will be abbreviated as
\begin{equation}
    \lagrangian_{zz} \equiv \DDt{\objective}{\des}. 
\end{equation}
In the case of inexact Newton-Krylov methods, where only the Hessian-vector products
\begin{equation}
\begin{split}
    %\DDt{\objective}{\des} \mbf{v}
    \lagrangian_{zz}
    &=
    %\left( \objective_{\des \des} + \dual^T \residuals_{\des \des} \right)
    \lagrangian_{\des \des}
    \mbf{v}
    + %\left ( 
      %  \objective_{\des \state} + \dual^T \residuals_{\des \state} 
      %\right) 
      \lagrangian_{\des \state}
      \impstate_\des
    \mbf{v}
    + \impstate_\des^T
      \lagrangian_{\state \des}
      %\left( \objective_{\state \des}  + \dual^T \residuals_{\state \des}
      %\right)
    \mbf{v}
    + \impstate_\des^T
      %\left( \objective_{\state \state} + \dual^T \residuals_{\state \state} \right)
      \lagrangian_{\state \state}
      \impstate_\des
    \mbf{v}
\end{split}
\end{equation}
%are required, we can apply the Hessian at the cost of two linear system solves as we apply $\residuals_\state^{-1}$ and $\residuals_\state^{-T}$. 
are required, we can apply the Hessian at the cost of two linear system solves.
The first solve forms
\begin{equation}
    \mbf{r}
    \coloneqq \impstate_\des \mbf{v} = \residuals_\state^{-1} \left(\residuals_\des \mbf{v} \right),
\label{eq:hess_vec_1}
\end{equation}
followed by 
\begin{equation}
\begin{split}
    \mbf{s}
    &\coloneqq \impstate_\des^T
      %\left( \objective_{\state \des}  + \dual^T \residuals_{\state \des} \right)
      \lagrangian_{\state \des}
    \mbf{v}
    + \impstate_\des^T
      %\left( \objective_{\state \state} + \dual^T \residuals_{\state \state} \right)
      \lagrangian_{\state \state}
      \impstate_\des
    \mbf{v}
    \\&=
    \residuals_{\state}^{-T}
    \left(
    \residuals_{\des}^{T}
      %\left( \objective_{\state \des}  + \dual^T \residuals_{\state \des}      \right)
      \lagrangian_{\state \des}
    \mbf{v}
    + 
    \residuals_{\des}^{T}
      %\left( \objective_{\state \state} + \dual^T \residuals_{\state \state} \right)
      \lagrangian_{\state \state}
    \mbf{r}
    \right),
\end{split}
\label{eq:hess_vec_2}
\end{equation}
to assemble into
\begin{equation}
\begin{split}
    \lagrangian_{zz} \mbf{v}
    &=
    %\left( \objective_{\des \des} + \dual^T \residuals_{\des \des} \right)
    \lagrangian_{\des \des}
    \mbf{v}
    + %\left ( 
      %  \objective_{\des \state} + \dual^T \residuals_{\des \state} 
      %\right)
      \lagrangian_{\des \state}
    \mbf{r}
    + \mbf{s}
    .
\end{split}
\label{eq:hess_vec_final}
\end{equation}
Hicken \cite{Hicken2014} provides tolerance bounds on those two solves such that superlinear convergence is observed for the optimization problem.
%as a two-step process where we first solve for $\mbs{\gamma} = \impstate_\des \mbf{v}$ using Eq. \eqref{eq:first_order_flow_sensitivity}
%\begin{equation}
%    \mbs{\gamma} = \impstate_\des \mbf{v} = - \residuals_\state^{-1} \residuals_\des \mbf{v}
%\end{equation}
%and then for $\mbs{\lambda} = \impstate_\des \mbf{v}$

\subsubsection{Reduced space implementation}
The state variables $\state$ are obtained through backward Euler fixed-point iterations described in Eq. \eqref{eq:backward_euler}, followed by the adjoint linear solution to Eq. \eqref{eq:adjoint_equation} before every control search direction computation.

The Newton search direction $\search{\des}$ is obtained through
\begin{equation}
    \lagrangian_{zz} 
    %\DDt{\objective}{\des}
    \search{\des} =
    - \mbf{g}_z 
    %\Dt{\objective}{\des}
    .
\end{equation}
%Note that the state variables $\state$ are obtained through some fixed-point iterations provided by the flow solver, followed by the adjoint linear solution before every search computation.
GMRES is used to solve the system by forming the matrix-vector products described in Eqs.~(\ref{eq:hess_vec_1}-\ref{eq:hess_vec_final}) with an identity-initialized BFGS preconditioner.
A relative tolerance of \oneEm{6} is used to determine convergence of the linear system.
In the case of the quasi-Newton implementation, the reduced-Hessian is completely replaced with the BFGS approximation.
%A line-search also complements the algorithm for globalization purposes.
A backtracking line-search also complements the algorithm for globalization purposes.

%% file: 3c_precond_full_space.tex
\section{Full-space preconditioning}
\label{sec:full_space_precond}
The convergence of Krylov solvers are directly related to the conditioning of the system
whose convergence highly depends on the quality of the preconditioner.
Biros and Ghattas \cite{Biros2005a} proposed preconditioners based on the reduced-space factorization.
Their $\PFour$ and $\PTwo$ preconditioners are defined as
\begin{equation}
    \PFour = 
    \begin{bmatrix}
    \lagrangian_{\state \state} \residuals_{\state}^{-1}
    & \bzero & \identity_{\nstate}
    \\
    \lagrangian_{\des \state} \residuals_{\state}^{-1}
    & \identity_{\ndes} & \residuals_{\des}^T \residuals_{\state}^{-T}
    \\
    \identity_{\nstate} & \bzero & \bzero 
    \end{bmatrix}
    \begin{bmatrix}
    \residuals_{\state} & \residuals_{\des} & \bzero
    \\
    \bzero & \approxHessian &  \bzero
    \\
    \bzero & \bzero & \residuals_{\state}^T
    \end{bmatrix}
\label{eq:p4_preconditioner}
\end{equation}
\begin{equation}
    \PTwo = 
    \begin{bmatrix}
    \bzero & \bzero & \identity_{\nstate}
    \\
    \bzero & \identity_{\ndes} & \residuals_{\des}^T \residuals_{\state}^{-T}
    \\
    \identity_{\nstate} & \bzero & \bzero 
    \end{bmatrix}
    \begin{bmatrix}
    \residuals_{\state} & \residuals_{\des} & \bzero
    \\
    \bzero & \approxHessian &  \bzero
    \\
    \bzero & \bzero & \residuals_{\state}^T
    \end{bmatrix}
    =
    \begin{bmatrix}
    \bzero & \bzero & \residuals_\state^T
    \\
    \bzero & \approxHessian & \residuals_\des^T
    \\
    \residuals_\state & \residuals_\des & \bzero
    \end{bmatrix}
\label{eq:p2_preconditioner}
\end{equation}
where $\approxHessian$ refers to an approximation of the reduced Hessian.
The original Broyden-Fletcher-Goldfarb-Shanno (BFGS) quasi-Newton method  \cite{Broyden1969,Fletcher1970,Goldfarb1970,Shanno1970} will be used to update the preconditioner in this work.
Note that the original version is used and not limited-memory BFGS~\cite{Nocedal1980} since it offers better convergence properties and the storage relative to the number of control variables is tractable in ASO.

The eigenvalue spectrum of $\PFour^{-1}$ and $\PTwo^{-1}$ applied onto the KKT matrix $\KKTmat$ from Eq. \eqref{eq:kkt_system} will be the set of eigenvalues \footnote{Despite the mistakes in the original derivation \cite{Biros2005a} addressed in \ref{app:preconditioner_inverses}}
\begin{equation}
    \mathcal{E}(\PFour \KKTmat) = \mathcal{E}(\PTwo \KKTmat) 
    = 
    \left\{ 
    \left\{ \mathcal{E}(\approxHessian^{-1} \lagrangian_{zz}) \right\}
    ,  \left\{ \mathbf{\identity_s} \right\}
    ,  \left\{ \mathbf{\identity_s} \right\}
    \right\}
    .
\end{equation}

The application of $\PFour^{-1}$ requires 4 linear solves (two applications of $\residuals_{\state}^{-1}$ and two applications of $\residuals_{\state}^{-T}$)
, while the application of $\PTwo^{-1}$ only requires 2 linear solves while keeping the same system conditioning.
By replacing the application of $\residuals_{\state}^{-1}$ and $\residuals_{\state}^{-T}$ by their respective preconditioners $\Tilde{\residuals}_{\state}^{-1}$ and $\Tilde{\residuals}_{\state}^{-T}$, we obtain the 
$\PFourA^{-1}$ and $\PTwoA^{-1}$ preconditioners.

The Trilinos preconditioner library \texttt{Ifpack} \cite{ifpack-website} is used to form the flow Jacobian and adjoint flow Jacobian preconditioners.
The preconditioner is an additive Schwarz domain decomposition method with 1 overlap between each subdomain and a threshold-based incomplete lower-upper factorization is used for each subdomain matrix with a fill ratio of 2, an absolute threshold of $\oneEm{4}$, a relative threshold of 1.0, and a drop tolerance of $\oneEm{10}$.

%% file: 4_flow_constraints.tex
\section{Flow constraints}
The two-dimensional steady-state compressible Euler equations are used to model the flow and its PDEs in Cartesian coordinates $x_1$ and $x_2$ are described, using Einstein notation, by
\begin{equation} \label{eq:GE}
\residuals(\state,\bx(\des))
= \del{\mathbf{f}_i(\state,\bx(\des))}{x_i} = \bzero,
\end{equation} 
where the conservative state vector $\state$, inviscid flux vector $\mathbf{f}_i$, are defined as,
% non preprint version
%\begin{equation}
%\state = \begin{bmatrix} \rho & \rho v_1 & \rho v_2 & \rho E \end{bmatrix}^T
%\quad \text{and} \quad
%\mathbf{f}_i = 
%\begin{bmatrix} 
%\rho v_i&
%\rho v_i v_1 + p \delta_{i1} & \rho v_i v_2 + p \delta_{i2}   &
%\rho E v_i +p v_i
%\end{bmatrix}^T
%.
%\end{equation}
%% preprint version
\begin{equation}
\state = \begin{bmatrix} \rho & \rho v_1 & \rho v_2 & \rho E \end{bmatrix}^T
\end{equation}
and
\begin{equation}
\mathbf{f}_i = 
\begin{bmatrix} 
\rho v_i&
\rho v_i v_1 + p \delta_{i1} & \rho v_i v_2 + p \delta_{i2}   &
\rho E v_i +p v_i
\end{bmatrix}^T
.
\end{equation}
The density, velocities, Kronecker delta function, and total energy are respectively denoted as $\rho$,  $v_i$, $\delta _{ij}$, and $E$.
The total energy is given by $E=e+\frac{1}{2}(\|\mbf{v}\|^2)$. The pressure $p$ is determined by the equation of state
\begin{equation}
p=(\gamma -1)\rho \left( E - \|\mbf{v}\|^2/2 \right) ,
\end{equation}
where $\gamma$ is the ratio of specific heats. 

The equations are discretized using the weak form of the discontinuous Galerkin method \cite{Reed1973} with a Roe flux \cite{Roe1981}.
The discretized equations are solved numerically using the Parallel High-Order Library for PDEs (\texttt{PHiLiP}) code which uses \texttt{deal.II} \cite{Bangerth2007} as its backbone. Trilinos' \texttt{Sacado} package \cite{sacado-website} is used to automatically differentiate the discretized residual to obtain the first and second partial derivatives 
$
\residuals_\state, 
\residuals_\bx, 
\dual^T \residuals_{\state \state}, 
\dual^T \residuals_{\state \bx},$ 
and $ \dual^T \residuals_{\des \bx}$.

%% file: 5_design_parametrization.tex
\section{Design parametrization}
While early works used the surface mesh points as control variables \cite{Jameson1988,Jameson1995}, increasingly complex geometric constraints arose and more intuitive controls were needed.
A popular method to parametrize the surface is to use free-form deformation (FFD) \cite{Sederberg1986}.
The displacements of the FFD control points are basically interpolated to locations within the FFD control box using B\'ezier curves.
The surface of interest is surrounded by a FFD box defined by an origin $\mbf{Q}$, and three vectors defining the parallepiped $\mbf{S}$, $\mbf{T}$, and $\mbf{U}$.
Any initial grid coordinate $\mbf{x}_{\text{init}} = (x_1, x_2, x_3)$ within the box has a corresponding local $(s,t,u)$ coordinate such that
\begin{equation}
\mbf{x}_{\text{init}} = \mbf{Q} + s \mbf{S} + t \mbf{T} + u \mbf{U}.
\end{equation}
The local coordinates can be obtained using some linear algebra
\begin{equation}
s = \frac
    {\mbf{T}\times\mbf{U}(\mbf{x}_{\text{init}}-\mbf{Q})}
    {\mbf{T}\times\mbf{U} \cdot \mbf{S}}
,\quad
t = \frac
    {\mbf{S}\times\mbf{U}(\mbf{x}_{\text{init}}-\mbf{Q})}
    {\mbf{S}\times\mbf{U} \cdot \mbf{T}}
\quad \text{and} \quad
s = \frac
    {\mbf{S}\times\mbf{T}(\mbf{x}_{\text{init}}-\mbf{Q})}
    {\mbf{S}\times\mbf{T} \cdot \mbf{U}}
.
\end{equation}
In two dimensions, it simplifies to
\begin{equation}
\mbf{x}_{\text{init}} = \mbf{Q} + s \mbf{S} + t \mbf{T}
,\quad
s = \frac
    {(\mbf{x}_{\text{init}} - \mbf{Q}) \cdot \mbf{T}^{\perp}}
    {\mbf{S} \cdot \mbf{T}^{\perp}}
\quad \text{and} \quad
t = \frac
    {(\mbf{x}_{\text{init}} - \mbf{Q}) \cdot \mbf{S}^{\perp}}
    {\mbf{T} \cdot \mbf{S}^{\perp}}
.
\end{equation}
where the vector $\mbf{v} = [v_1, v_2]$ has its corresponding perpendicular vector $\mbf{v}^\perp = [v_2, -v_1]$.
Control points $\mbf{P}_{ijk}$ of the FFD parametrization are evenly spaced on the local $(s,t,u)$ coordinate such that for $(n_s, n_t, n_u)$ control points in their respective directions
\begin{equation}
    \mbf{P}_{ijk} = \mbf{Q}
    + \frac{i}{n_s} \mbs{S}
    + \frac{j}{n_t} \mbs{S}
    + \frac{k}{n_u} \mbs{S}
    .
\end{equation}
Finally,
we can move the control points $\mbf{P}_{ijk}$ and recover the new location of a point $(x,y,z)$ with corresponding local coordinates $(s,t,u)$
through the deformation function defined by the tensor product Bernstein polynomial
\begin{equation}
    \mbf{x}_{surf} = 
    \sum_{i=0}^{n_s} \binom{n_s}{i}(1-s)^{n_s - i} s^{i}
    \left[
    \sum_{j=0}^{n_t} \binom{n_t}{j}(1-t)^{n_t - j} t^{j}
    \left[
    \sum_{k=0}^{n_u} \binom{n_u}{k}(1-u)^{n_u - k} u^{k}
    \mbf{P}_{ijk}
    \right]
    \right]
    .
\end{equation}
A subset of the control points $\mbf{P}_{ijk}$ will be chosen as the control variables $\des$.
Moreover, the necessary derivatives of the surface mesh points with respect to the control points are simply given by
\begin{equation}
    \left(\del{ \mbf{x}_{surf}}{ \mbf{P}_{ijk}}\right)_{ab}
    =
    \binom{n_s}{i}(1-s)^{n_s - i} s^{i}
    \binom{n_t}{j}(1-t)^{n_t - j} t^{j}
    \binom{n_u}{k}(1-u)^{n_u - k} u^{k}
    .%\delta_{ab},
\end{equation}
%where $\delta_{ab}$ is the Kronecker delta function.

%The FFD parametrization maps the initial grid coordinates $(x,y,z)$ onto its own local $(s,t,u)$ coordinate system based on the bounding FFD box.
%The FFD box control
%Since the FFD parametrization maps the initial grid onto its own coordinate system

%% file: 6_meshmover.tex
\section{Volume mesh movement}
%Surface node displacements required by the previously defined mesh movement are recovered from the displacements of the FFD control points.
%The FFD control points are used as control variables.
Once the surface displacements have been recovered through the FFD, the volume mesh nodes need to be displaced to ensure a valid mesh.
The high-order mesh is deformed using a linear elasticity mesh movement \cite{Tezduyar1992,Truong2008,Persson2009,Brown2018} that propagates prescribed surface deformations into the domain.

The elastic equations with Dirichlet boundary conditions are described by
\begin{subequations}
\begin{equation}
    -\del{}{x_i} \lambda \left( \del{(\Delta \xvol)_j}{x_j} \right )
    - \del{}{x_j} \left (\mu \del{(\Delta \xvol)_j}{x_i} \right)
    - \del{}{x_j} \left( \mu  \del{(\Delta \xvol)_i}{x_j}  \right)
    = b_i
    ,\quad \bx \in \domain
    ,
\end{equation}
\begin{equation}
%\quad \text{and} \quad
     \left(\Delta \xvol \right)_i = \mbf{M}_{ij}\left(\Delta \xsurf^{\text{BC}}\right)_j
    ,\quad \bx \in \boundary
    ,
\label{eq:linear_elasticity}
\end{equation}
\end{subequations}
where $(\Delta \xvol)$ represents the volume displacements, $(\Delta \xsurf)$ represent the surface displacements, and $\mbf{M}_{ij}$ maps the surface indices to their respective volume degree of freedom indices. Additionally, $\lambda$ and $\mu$ are the first and second Lam\'e parameters, and $\mbf{b}$ are the external forces.
The equation is discretized using the continuous Galerkin method since the node locations must be continuous across the elements.
%Once Eq. \eqref{eq:linear_elasticity} is discretized, we can form the resulting system
%\begin{equation}
%    \mbf{S}\Bar{\mbf{y}} = \bzero \quad \text{and} \quad \mbf{y} = \mbf{C} \Bar{\mbf{y}} + \mbf{k}
%\end{equation}
%where $\mbf{S}$ is the stiffness matrix with unconstrained degrees of freedom such that no boundary conditions have been applied and $\Bar{\mbf{y}}$ is the solution to the unconstrained system.
%The Dirichlet boundary conditions and any other affine constraints coming from hanging nodes are represented by the matrix of homogeneous coefficients $\mbf{C}$ and the vector of inhomogeneous values $\mbf{k}$.
%
%The matrix $\mbf{C}$ contains the homogeneous part of the affine constraints
Since we only have boundary displacements, there are no external forces, such that $\mbf{b} = \bzero$.

The Dirichlet boundary conditions are imposed by zero-ing the corresponding row of the stiffness matrix, setting the diagonal to 1, and the right-hand side to the Dirichlet value.
We emphasize the importance of the described Dirichlet boundary conditions implementation since the resulting system to be solved becomes
\begin{equation}
    \mbf{S} \left(\Delta \xvol\right) = \mbf{M}\left(\Delta \xsurf^{\text{BC}}\right)
    ,
\end{equation}
where $\mbf{S}$ represents the stiffness matrix resulting from the discretization.
Therefore, the mesh sensitivities are simply given by
\begin{equation}
    \del{\xvol}{\xsurf} = \mbf{S}^{-1} \mbf{M}
    .
\end{equation}

Since the mesh sensitivities require a linear solve every time either $\objective_{\des}$,
$\objective_{\des \state}$,
$\objective_{\state \des}$,
$\objective_{\des \des}$,
$\residuals_{\des}$,
$\residuals_{\des \state}$,
$\residuals_{\state \des}$,
or 
$\residuals_{\des \des}$
is applied onto a vector, storing the chain-ruled derivative at the cost of $\ndes$ linear solves
\begin{equation}
    \del{\xvol}{\des} 
    = \del{\xvol}{\xsurf}\del{\xsurf}{\des} 
    =  \mbf{S}^{-1} \mbf{M} \del{\xsurf}{\des}
\end{equation}
will result in higher performance since the total number of the mesh sensitivity applications will likely exceed the number of control variables.
Furthermore, since the control variables and mesh movement always stem from an original mesh, the chain-ruled sensitivities only have to be pre-processed once.
However, the memory requirement will be equivalent to storing $\ndes$ meshes.
%While it may be tempting to store the chain-rule derivative
%\begin{equation}
%    \del{\xvol}{\des} = \del{\xvol}{\xsurf}\del{\xsurf}{\des} =  \mbf{S}^{-1} \mbf{M_{ij}} (\del{xsurf}{\des}
%\end{equation}

%and can easily be applied onto a single vector.         
%Once the equations have been discretized, we can assemble the eb

%A result of applying the Dirichlet this way is that 
%After setting the forces $\mbf{b}$ to zero and assigning Dirichet boundary conditions coming from the surface displacements, solving Eq. \eqref{eq:linear_elasticity} recovers the volume displacements of the high-order nodes.

%% file: 7_cost_breakdown.tex
\section{Cost analysis}
\newcommand{\work}[1]{W\left( #1 \right)}
\newcommand{\workR}[1]{W_{\residuals}\left( #1 \right)}
\newcommand{\vmult}[1]{\texttt{vmult}\left(#1\right)}
%The assembly of a single flow residual will be used as a unit of work.
%Its derivatives through automatic differentiation (AD) have been analyzed as a relative cost to assemble the actual residual \cite{AndreasGriewank2008}.
%We will see that depending on the framework used, it may be advantageous to form the sensitivity matrices explicitly Depending on multiple factors

\newcommand{\nb}[0]{n_p}
\subsection{Residual cost}

The analysis is done in $d$-dimensions for the weak form of DG discretizing the $s=d+2$ Euler equations. 
We use an arbitrary tensor-product basis of order $p$ for the solution, and a Gauss-Lobatto tensor-product basis of order $p$ for the grid.
A total of $\nb = (p+1)^d$ basis are used to discretize the solution and mesh.
The integral computations use $\nb$ cubature nodes and $n_f = (p+1)^{d-1}$ quadrature nodes.
Although the solution and the grid are discretized through tensor-product bases, the current implementation does not use the more cost effective sum-factorization method \cite{Orszag1980_SumFact}. 

The number of flops required to evaluate the domain contribution is around 
$(32\nb^2 + 45\nb)$
in 2D 
and 
$(58\nb^2 + 121\nb)$,
in 3D 
as per \ref{app:domain_residual_flops}, 
whereas the surface contribution is around 
$(32\nb n_f + 227 n_f)$
in 2D 
and 
$(48 \nb n_f + 299 n_f)$
in 3D 
as detailed in  \ref{app:surface_residual_flops}.
For a structured 1/2/3D mesh, there will be 1/2/3 times as many faces as there are elements.
Therefore, we obtain a total residual cost of
\begin{equation}
\work{\residuals} = 32\nb^2 + 45\nb + 2(32\nb n_f + 227 n_f)
\end{equation}
or
\begin{equation}
\work{\residuals} = 32(p+1)^4 + 64 (p+1)^{3} + 45(p+1)^2 + 454 (p+1)
\end{equation}
in 2D and
\begin{equation}
\work{\residuals} = 58\nb^2 + 121\nb+ 3(48 \nb n_f + 299 n_f)
\end{equation}
or
\begin{equation}
\work{\residuals} = 58(p+1)^{6} + 144 (p+1)^{5} + 121 (p+1)^{3} + 897 (p+1)^{2}
\end{equation}
flops in 3D, where $\work{\,}$ denotes the total number of flops.
The absolute flops count is tabulated in Table \ref{tab:total_residual_cost}.
Additionally, $\workR{\,}$ will be used to denote other type of operations relative to $\residuals$ such that $\workR{\texttt{ops}} = \work{\texttt{ops}}/\work{\residuals}$ and $\workR{\residuals} = 1$.

\input{table/cost_residual_table}

\subsection{Matrix-vector products}

The flops necessary to evaluate the matrix-vector products of the Jacobian or the Hessians are directly proportional to the sparsity of those operators.
However, the time required to effectuate those flops differ from the flops of the residual assembly.
Through numerical tests, we determined that the flops associated with the matrix-vector products are approximately three times faster than the ones from the residual assembly.
The incongruity between the number of flops and timings can have multiple causes such as memory accesses, cache locality, and SIMD optimizations, which will vary with the implementation of the residual assembly.
Therefore, the absolute work $\work{}$ and relative work $\workR{}$ in Table 
\ref{tab:dRdW_vmult_cost}, 
\ref{tab:dRdX_vmult_cost}, 
and \ref{tab:d2RdXdX_vmult_cost}
will be divided by 3 in the remainder of the paper to more accurately reflect the timings.
The experimental timings used a single core of a Ryzen 2700X processor with DDR4-2933 MHz memory to compare the residual assembly with the sparse matrix-vector products within \texttt{PHiLiP}.

\subsubsection{Jacobian-vector product}
A matrix-vector product of an already formed matrix will only involve the non-zeros of the matrix.
The number of nonzeros per row of $\residuals_{\state}$ is equal to the cell stencil size times the number of degrees of freedom per cell
\begin{equation}
    N_{nz}^{\residuals_{\state}}(p,d) = (1 + 2d) (d+2) (p+1)^d.
\end{equation}
Therefore, the residuals corresponding to one cell has $\nb s$ rows and $N^{\residuals_\state}_{nz}(p,d)$ non-zero columns.
The total flops required to form the block-row-vector product is given by
\begin{equation}
\begin{split}
    \work{\vmult{\residuals_{\state} \mbf{v}}}
    &= 2 \nb s N^{\residuals_\state}_{nz}(p,d) 
    \\&= 2 \nb s (1 + 2d) (d+2) (p+1)^d 
    %\\&= 2 (1 + 2d) (d+2)^2 (p+1)^{2d}
    \\&= (4d^3+18d^2+24d+8) (p+1)^{2d}.
\end{split}
\end{equation}
\input{table/cost_vmult_table}

A similar analysis can be used to obtain $\work{\vmult{\residuals_{\xvol} \mbf{v}}}$.
This time, the cell's corresponding residuals only depends on its own nodes.
The resulting work is therefore
\begin{equation}
\begin{split}
    \work{\vmult{\residuals_{\xvol} \mbf{v}}} 
    &= 2 \nb s d (p+1)^d 
    = (2d^2 + 4d) (p+1)^{2d}
.
\end{split}
\end{equation}
The total and relative work $\work{}$ and $\workR{}$ are tabulated in Table \ref{tab:dRdW_vmult_cost} and \ref{tab:dRdX_vmult_cost}, which will be divided by 3 as previously explained.
Furthermore, the flow and adjoint Jacobian preconditioner use an ILUT fill ratio of 2 (not to be confused with ILU(2)), it has twice as many nonzeros, and its application is therefore twice as expensive.
% \red{
% Although the flops required to perform the matrix-vector product $\vmult{\residuals_\state \mbf{v}}$ is larger than the flops required to assemble the residual, our numerical experiments have shown that the $\vmult{}$ wallclock timings are approximately a third of the estimated costs.
% The incongruity between the number of flops and timings can have multiple causes such as memory accesses, cache locality, and SIMD optimizations.
% The total cost calculations used later will use the more representative third of the cost of the ones shown in Table \ref{tab:dRdW_vmult_cost} and \ref{tab:dRdX_vmult_cost}.
% }
%likely due to the different nature of flops being performed.
%the highly efficient matrix-vector product as opposed to 
%This is likely due to the memory locality and parallelization

\subsubsection{Hessian-vector product}
The dual-weighted residual Hessian block $\dual^T \residuals_{\state \state}$ has the same sparsity pattern as $\residuals_{\state}$ and therefore has the same matrix-vector product cost as described in Table \ref{tab:dRdW_vmult_cost}.
Similarly, the blocks $\dual^T \residuals_{\state \xvol}$ and $\dual^T \residuals_{\xvol \state}$ have the same sparsity pattern as $\residuals_{\xvol}$ and the same matrix-vector product cost described in Table \ref{tab:dRdX_vmult_cost}.
Finally, the matrix-vector product cost of $\dual^T \residuals_{\xvol \xvol}$ listed in Table \ref{tab:d2RdXdX_vmult_cost} is $d/(d+2)$ times the matrix-vector product cost of $\dual^T \residuals_{\state \xvol}$ since it  contains less rows.
%\red{Similarly to the previous section, }

\begin{table}[h]
\centering
\begin{tabular}{
    l 
    S[table-format=7.0] @{\,\,\( / \)\,} S[table-format=3.2]
    |
    S[table-format=7.0] @{\,\,\( / \)\,} S[table-format=3.2]
    }
  & \multicolumn{4}{c}{d}\\
\cmidrule{2-5}
p & \multicolumn{2}{c}{2} & \multicolumn{2}{c}{3} \\
\midrule
1 & 128 & 0.06 & 1152 & 0.09 \\ 
2 & 648 & 0.11 & 13122 & 0.15 \\ 
3 & 2048 & 0.14 & 73728 & 0.18 \\ 
4 & 5000 & 0.16 & 281250 & 0.20 \\ 
\end{tabular}
\caption{
Flops/Relative work required to assemble matrix-vector product $\vmult{\residuals_{\xvol \xvol} \mbf{v}}$ when $\residuals_{\xvol \xvol}$ has been assembled already.
$\work{\vmult{\residuals_{\xvol \xvol} \mbf{v}}} / \workR{\vmult{\residuals_{\xvol \xvol} \mbf{v}}}$ 
}
\label{tab:d2RdXdX_vmult_cost}
\end{table}

\subsection{Automatic differentiation}
\subsubsection{Forward and reverse work}
\newcommand{\relativework}[0]{\omega}
\newcommand{\seed}[0]{\mbf{S}}
Given the work required to assemble the residual $\work{\residuals}$ and a seed matrix $\seed \in \real^{q \times k}$, where $q$ represents the number of independent variables and $k$ represents the number of columns on which to apply the residual derivative, a time complexity analysis of automatic differentiation \cite{AndreasGriewank2008} shows that the relative work $\workR{\residuals_{\state} \seed}$ of the vector forward mode
is given by
\begin{equation}
1+ k \leq \workR{\residuals_{\state} \seed} \leq 1.5 + k
,
\label{eq:forward_relwork}
\end{equation}
where the range depends on the cost of memory access.
The lower relative complexity occurs when the residual assembly is memory-bound, whereas the upper bound occurs when the residual assembly is compute-bound.
For the remainder of this paper, we use an average $\workR{\residuals_{\state} \seed}\approx 1.25+k$.
The vector reverse mode's relative work $\workR{\seed^T \residuals_{\state}}$
is given by
\begin{equation}
1 + 2k \leq \workR{\seed^T \residuals_{\state}} \leq 1.5 + 2.5k
\label{eq:gradk_relwork}
\end{equation}
as derived in \cite{AndreasGriewank2008}.
For the remainder of this paper, we use an average $\workR{\seed^T \residuals_{\state}} \approx 1.25+2.25k$.

\subsubsection{Matrix-free product}
If only a Jacobian-vector product is needed, then $k=1$ and 
\begin{equation}
    \workR{ \residuals_{\state} \seed } = 2.25.
\end{equation}
The same applies to the mesh Jacobian since the AD relative cost only depends on the number of columns the derivative is applied onto such that 
\begin{equation}
    \workR{ \residuals_{\xvol} \seed } = 2.25.
\end{equation}
Similarly, the transposed versions with $k=1$ in the reverse mode results in
\begin{equation}
    \workR{ \seed^T \residuals_{\state} } =
    \workR{ \seed^T \residuals_{\xvol} } = 3.5.
\end{equation}

In the case of the dual-weighted residual Hessian, we simply use the reverse mode to first evaluate the dual-weighted Jacobian
\begin{equation}
    \workR{ \dual^T \residuals_{\state} } = 
    \workR{ \dual^T \residuals_{\xvol} } = 3.5.
    \label{eq:reverse_dual_residual_cost}
\end{equation}
The forward mode is then applied onto the seed vector with $k=1$ to multiply the above cost by 2.25.
% non preprint version
%\begin{equation}
%    \workR{ \dual^T \residuals_{\state \state} \seed } = 
%    \workR{ \dual^T \residuals_{\state \xvol} \seed } = 
%    \workR{ \dual^T \residuals_{\xvol \state} \seed } = 
%    \workR{ \dual^T \residuals_{\xvol \xvol} \seed } = 7.875.
%    \label{eq:reverse_forward_dual_residual_cost}
%\end{equation}
% preprint version
\begin{equation}
\begin{split}
    \workR{ \dual^T \residuals_{\state \state} \seed }
    &=
    \workR{ \dual^T \residuals_{\state \xvol} \seed } 
    \\&=
    \workR{ \dual^T \residuals_{\xvol \state} \seed } 
    \\&=
    \workR{ \dual^T \residuals_{\xvol \xvol} \seed }
    \\&=
    7.875.
    \label{eq:reverse_forward_dual_residual_cost}
\end{split}
\end{equation}

\subsection{Forming the derivatives}
While matrix-free preconditioners for DG are in development \cite{Pazner2018,Franciolini2020}, more standard algorithms, such as ILU, require the explicit matrix to be available.
To build those derivatives, we simply have to apply the forward mode onto a seed matrix $\seed$ with columns formed by the standard bases corresponding to independent variables.

In the case of the volume work, the number of independent variables are the number of degrees of freedom within the cell $k=(d+2)(p+1)^d$.
For the face computations, the two neighbouring cells' degrees of freedom are the independent variables, totalling $k=2(d+2)(p+1)^d$ columns of $\seed$.
The relative work $\workR{\residuals_{\state}}$ therefore changes with the dimension and polynomial degree and is given by Table \ref{tab:dRdW_cost}.
A similar analysis can be performed on $\workR{\residuals_{\xvol}}$ by only using its own nodes to count the number of independent variables to obtain Table \ref{tab:dRdX_cost}.
\begin{table}[h]
\begin{minipage}{0.48\textwidth}
\centering
    \begin{tabular}{l 
    S[table-format=4.1]
    |
    S[table-format=4.1]
    }
      & \multicolumn{2}{c}{d}\\
    \cmidrule{2-3}
p & \multicolumn{1}{c}{2} & \multicolumn{1}{c}{3} \\
\midrule
1 & 28.0 & 66.7 \\ 
2 & 55.5 & 201.9 \\ 
3 & 90.8 & 448.4 \\ 
4 & 134.0 & 838.1 \\ 
    \end{tabular}
    \caption{Relative work to assemble the flow Jacobian explicitly $\workR{\residuals_{\state}}$ using AD.}
    \label{tab:dRdW_cost}
\end{minipage}
 \hfill
\begin{minipage}{0.48\textwidth}
\centering
    \begin{tabular}{l 
    S[table-format=4.1]
    |
    S[table-format=4.1]
    }
      & \multicolumn{2}{c}{d}\\
    \cmidrule{2-3}
p & \multicolumn{1}{c}{2} & \multicolumn{1}{c}{3} \\
\midrule
1 &  9.2 & 25.2 \\ 
2 & 19.2 & 82.2 \\ 
3 & 33.2 & 193.2 \\ 
4 & 51.2 & 376.2 \\ 
    \end{tabular}
    \caption{Relative work to assemble the mesh Jacobian explicitly $\workR{\residuals_{\xvol}}$ using AD.}
    \label{tab:dRdX_cost}
\end{minipage}
\end{table}

When it comes to finding the cost to form the dual-weighted residual Hessian blocks, we can start with the cost from Eq. \eqref{eq:reverse_dual_residual_cost}.
The seed matrix in Eq. \eqref{eq:reverse_forward_dual_residual_cost} is then the same as the ones described to form $\residuals_{\state}$ and $\residuals_{\xvol}$.
As a result, the costs of assembling the dual-weighted residual Hessian blocks are
\begin{equation}
    \workR{\dual^T \residuals_{\state \state}} = 
    \workR{ \dual^T \residuals_{\state} }
    \times
    \workR{\residuals_{\state}}
    = 3.5 \workR{\residuals_{\state}}
,
\end{equation}
%non preprint
%\begin{equation}
%    \workR{\dual^T \residuals_{\state \xvol}} 
%    = \workR{\dual^T \residuals_{\xvol \xvol}}
%    = \workR{ \dual^T \residuals_{\xvol} }
%    \times
%    \workR{\residuals_{\xvol}}
%    = 3.5 \workR{\residuals_{\xvol}}
%.
%\end{equation}
%preprint split
\begin{equation}
\begin{split}
    \workR{\dual^T \residuals_{\state \xvol}} 
    &= \workR{\dual^T \residuals_{\xvol \xvol}}
    \\&= \workR{ \dual^T \residuals_{\xvol} }
    \times
    \workR{\residuals_{\xvol}}
    \\&= 3.5 \workR{\residuals_{\xvol}}
.
\end{split}
\end{equation}
Note that $\workR{\dual^T \residuals_{\state \xvol}}$ can be assembled more cheaply than $\workR{\dual^T \residuals_{\xvol \state}}$ and they are simply a transpose of each other.

\begin{table}[h]
\begin{minipage}{0.48\textwidth}
\centering
    \begin{tabular}{l 
    S[table-format=4.1]
    |
    S[table-format=4.1]
    }
      & \multicolumn{2}{c}{d}\\
    \cmidrule{2-3}
p & \multicolumn{1}{c}{2} & \multicolumn{1}{c}{3} \\
\midrule
1 & 98.0 & 233.5 \\ 
2 & 194.3 & 706.5 \\ 
3 & 317.7 & 1569.5 \\ 
4 & 468.9 & 2933.3 \\ 
    \end{tabular}
    \caption{Relative work to assemble the dual-weighted residual Hessian explicitly with second order state variables derivatives $\workR{\dual^T \residuals_{\state \state}}$ using AD.}
    \label{tab:d2RdWdW_cost}
\end{minipage}
 \hfill
\begin{minipage}{0.48\textwidth}
\centering
    \begin{tabular}{l 
    S[table-format=4.1]
    |
    S[table-format=4.1]
    }
      & \multicolumn{2}{c}{d}\\
    \cmidrule{2-3}
p & \multicolumn{1}{c}{2} & \multicolumn{1}{c}{3} \\
\midrule
1 & 32.4 & 88.4 \\ 
2 & 67.4 & 287.9 \\ 
3 & 116.4 & 676.4 \\ 
4 & 179.4 & 1316.9 \\ 
    \end{tabular}
    \caption{Relative work to assemble the dual-weighted residual Hessian explicitly with a mesh derivative $\workR{\dual^T \residuals_{\state \xvol}}$ and $\workR{\dual^T \residuals_{\xvol \xvol}}$ using AD.}
    \label{tab:d2RdX_cost}
\end{minipage}
\end{table}

\subsection{Matrix-free versus forming the matrices}
\newcommand{\worthform}[1]{C\left( #1 \right)}
The matrix-free approaches to evaluate Jacobian-vector and transpose-Jacobian-vector products has a constant relative cost of 2.25 and 3.5 residual evaluations.
The matrix-free dual-weighted residual Hessian-vector products also have a relative constant complexity of 7.875.
Those matrix-free products are more expensive than performing the matrix-vector product when the matrix has been pre-assembled.
It is not obvious whether it is worth investing a large amount of time once to save flops by performing the matrix-vector products with a pre-assembled matrix, rather than forming the matrix-vector products on-the-fly through AD.
In the case of the flow Jacobian $\residuals_{\state}$, the explicit matrix is required to form the preconditioner of the forward and adjoint problem.

However, for $\residuals_{\xvol}$, $\dual^T \residuals_{\state \state}$, $\dual^T \residuals_{\state \xvol}$, $\dual^T \residuals_{\xvol \state}$, and $\dual^T \residuals_{\xvol \xvol}$, we determine how many times the matrices have to be used by dividing the initial cost by the cost savings
\begin{equation}
    \worthform{\residuals_{\xvol}} = 
    \frac{
    \workR{\residuals_{\xvol}} 
    }{
    \left (
        \workR{\vmult{\residuals_{\xvol} \mbf{v}}} 
        - \workR{\residuals_{\xvol} \mbf{v}} 
    \right)
    }
.
\end{equation}
The number of times the matrix-products have to be formed to recover the initial cost are listed in Tables \ref{tab:dRdX_form_worth}--\ref{tab:d2RdXdX_form_worth}.
The optimization results in the following section are in 2D up to $p=3$.
Since the matrix assembly cost is recouped after very few matrix-vector products and the total number of matrix-vector products required is not known \emph{a priori}, the derivative operators are always fully formed rather than assembling the matrix-vector product on-the-fly.

\begin{table}[h]
\begin{minipage}{0.48\textwidth}
\centering
    \begin{tabular}{l 
    S[table-format=4.0]
    |
    S[table-format=4.0]
    }
      & \multicolumn{2}{c}{d}\\
    \cmidrule{2-3}
p & \multicolumn{1}{c}{2} & \multicolumn{1}{c}{3} \\
\midrule
1 &    2 &    6 \\ 
2 &    4 &   19 \\ 
3 &    8 &   45 \\ 
4 &   12 &   88 \\ 
    \end{tabular}
    \caption{Number of times we need the matrix-vector product $\residuals_{\xvol} \mbf{v}$ to warrant forming the matrix explicitly instead of assembling the matrix-product on-the-fly with AD. Its assembly cost has been divided by 2 since once assembled, it can be used for the required transpose matrix-vector product.}
    \label{tab:dRdX_form_worth}
\end{minipage}
 \hfill
\begin{minipage}{0.48\textwidth}
\centering
    \begin{tabular}{l 
    S[table-format=4.0]
    |
    S[table-format=4.0]
    }
      & \multicolumn{2}{c}{d}\\
    \cmidrule{2-3}
p & \multicolumn{1}{c}{2} & \multicolumn{1}{c}{3} \\
\midrule
1 &   13 &   32 \\ 
2 &   27 &  102 \\ 
3 &   46 &  234 \\ 
4 &   69 &  447 \\ 
    \end{tabular}
    \caption{Number of times we need the matrix-vector product $\residuals_{\state \state} \mbf{v}$ to warrant forming the matrix explicitly instead of assembling the matrix-product on-the-fly with AD.}
    \label{tab:d2RdWdW_form_worth}
\end{minipage}
\end{table}

\begin{table}[h]
\begin{minipage}{0.48\textwidth}
\centering
    \begin{tabular}{l 
    S[table-format=4.0]
    |
    S[table-format=4.0]
    }
      & \multicolumn{2}{c}{d}\\
    \cmidrule{2-3}
p & \multicolumn{1}{c}{2} & \multicolumn{1}{c}{3} \\
\midrule
1 &    2 &    6 \\ 
2 &    4 &   18 \\ 
3 &    7 &   44 \\ 
4 &   12 &   85 \\ 
    \end{tabular}
    \caption{Number of times we need the matrix-vector product $\residuals_{\state \xvol} \mbf{v}$ to warrant forming the matrix explicitly instead of assembling the matrix-product on-the-fly with AD. Its assembly cost has been divided by 2 since once assembled, it can be used for the required transpose matrix-vector product.}
    \label{tab:d2RdWdX_form_worth}
\end{minipage}
 \hfill
\begin{minipage}{0.48\textwidth}
\centering
    \begin{tabular}{l 
    S[table-format=4.0]
    |
    S[table-format=4.0]
    }
      & \multicolumn{2}{c}{d}\\
    \cmidrule{2-3}
p & \multicolumn{1}{c}{2} & \multicolumn{1}{c}{3} \\
\midrule
1 &    4 &   11 \\ 
2 &    9 &   37 \\ 
3 &   15 &   87 \\ 
4 &   23 &  169 \\ 
    \end{tabular}
    \caption{Number of times we need the matrix-vector product $\residuals_{\xvol \xvol} \mbf{v}$ to warrant forming the matrix explicitly instead of assembling the matrix-product on-the-fly with AD.}
    \label{tab:d2RdXdX_form_worth}
\end{minipage}
\end{table}

\newcommand{\mult}[2]{\the\numexpr (#1)*(#2) \relax}
\newcommand{\precondcost}[0]{2}
\newcommand{\LinearItCost}[0]{\the\numexpr \precondcost+1}
\newcommand{\LinearSolveIt}[0]{120}
\newcommand{\LinearSolveWork}[0]{\mult{\LinearSolveIt}{\LinearItCost}} 
\newcommand{\nBackwardEulerIt}[0]{4}
\newcommand{\BFGSWork}[0]{\mult{\the\numexpr \nBackwardEulerIt+1}{\LinearSolveWork}}

%To simplify our estimates, we will average out some numbers obtained from the test case described in the results section \ref{sec:results}.
%Solving to machine accuracy a system involving the flow Jacobian or its transpose required around \LinearSolveIt{}, which is equivalent to \LinearSolveWork{} work units.

%Within the reduced space, the backward Euler method requires around \nBackwardEulerIt{} nonlinear iterations to converge, followed by an adjoint solve.
%Therefore, a reduced-space quasi-Newton framework will require around \BFGSWork{} + $30 N_{nz}^{p,d}$ work units per design cycle, whereas the reduced-space Newton cycle will depend on the number of conjugate gradient (CG) iterations.
%A single CG iteration will ne
%Meanwhile, the the

%Furthermore, the backward Euler method use after each design cycle in the reduced space required around 4 nonlinear iterations to converge, 
%In addition to the The reduced quasi-Newton method therefore requires 

%% file: table/cost_residual_table.tex
\begin{table}[h]
    \centering
    \begin{tabular}{
    l
    % d1
    S[table-format=7.0] @{\,\,\( / \)\,}
    S[table-format=7.0] @{\,\,\( / \)\,} 
    S[table-format=7.0]
    |
    % d2
    S[table-format=7.0] @{\,\,\( / \)\,}
    S[table-format=7.0] @{\,\,\( / \)\,} 
    S[table-format=7.0]
    }
     & \multicolumn{6}{c}{d}
    \\ \cmidrule{2-7}
    p            & \multicolumn{3}{c}{2} & \multicolumn{3}{c}{3} \\
    \midrule
1  & 692  &  710  &  2112   & 4680  &  2732  &  12876  \\ 
2  & 2997  &  1545  &  6087   & 45549  &  14355  &  88614  \\ 
3  & 8912  &  2956  &  14824   & 245312  &  53936  &  407120  \\ 
4  & 21125  &  5135  &  31395   & 921375  &  157475  &  1393800  \\ 
    \end{tabular}
    \caption{Flops required to assemble the volume/face/total residual of a single cell.}
    \label{tab:total_residual_cost}
\end{table}

%% file: table/cost_vmult_table.tex
\begin{table}[h]
\centering
\begin{minipage}{0.48\textwidth}
\resizebox{\textwidth}{!}{%
\begin{tabular}{
    l 
    S[table-format=8.0] @{\,\,\( / \)\,} S[table-format=3.1]
    |
    S[table-format=8.0] @{\,\,\( / \)\,} S[table-format=3.1]
    }
  & \multicolumn{4}{c}{d}\\
\cmidrule{2-5}
p & \multicolumn{2}{c}{2} & \multicolumn{2}{c}{3} \\
\midrule
1 & 2560 &  1.2 & 22400 &  1.7 \\ 
2 & 12960 &  2.1 & 255150 &  2.9 \\ 
3 & 40960 &  2.8 & 1433600 &  3.5 \\ 
4 & 100000 &  3.2 & 5468750 &  3.9 \\  
\end{tabular}
}
\caption{
Flops/Relative work required to assemble matrix-vector product $\vmult{\residuals_\state \mbf{v}}$ when $\residuals_\state$ has been assembled already.
$\work{\vmult{\residuals_\state \mbf{v}}} / \workR{\vmult{\residuals_\state \mbf{v}}}$
}
\label{tab:dRdW_vmult_cost}
\end{minipage}
\begin{minipage}{0.48\textwidth}
\resizebox{\textwidth}{!}{%
\begin{tabular}{
    l 
    S[table-format=7.0] @{\,\,\( / \)\,} S[table-format=3.2]
    |
    S[table-format=7.0] @{\,\,\( / \)\,} S[table-format=3.2]
    }
  & \multicolumn{4}{c}{d}\\
\cmidrule{2-5}
p & \multicolumn{2}{c}{2} & \multicolumn{2}{c}{3} \\
\midrule
1 & 256 & 0.12 & 1920 & 0.15 \\ 
2 & 1296 & 0.21 & 21870 & 0.25 \\ 
3 & 4096 & 0.28 & 122880 & 0.30 \\ 
4 & 10000 & 0.32 & 468750 & 0.34 \\ 
\end{tabular}
}
\caption{
Flops/Relative work required to assemble matrix-vector product $\vmult{\residuals_\xvol \mbf{v}}$ when $\residuals_\xvol$ has been assembled already.
$\work{\vmult{\residuals_\xvol \mbf{v}}} / \workR{\vmult{\residuals_\xvol \mbf{v}}}$
}
\label{tab:dRdX_vmult_cost}
\end{minipage}
\end{table}

%\begin{table}[h]
%\centering
%    \begin{tabular}{l 
%    S[table-format=7.0]
%    |
%    S[table-format=7.0]
%    }
%      & \multicolumn{2}{c}{d}\\
%    \cmidrule{2-3}
%p & \multicolumn{1}{c}{1} & \multicolumn{1}{c}{2} \\
%\midrule
%1 & 2560 & 22400 \\ 
%2 & 12960 & 255150 \\ 
%3 & 40960 & 1433600 \\ 
%4 & 100000 & 5468750 \\ 
%    \end{tabular}
%    \caption{Flops/Relative Work $\work{\vmult{\residuals_\state \mbf{v}}} / \workR{\vmult{\residuals_\state \mbf{v}}}$ required to assemble matrix-vector product $\vmult{\residuals_\state \mbf{v}}$ when $\residuals_\state$ has been assembled already.}
%    \label{tab:vmult_cost}
%\end{table}

%% file: 7_results.tex
\section{Results}
\label{sec:results}
\subsection{Test case description}

\begin{figure}[h]
    \centering
    \subfloat[][Initial]{ \label{fig:initial_gaussian_bump}
        %\frame{
        \includegraphics
        [width=0.8\textwidth,trim={0.5cm 8.0cm 10.5cm 6cm},clip]
        {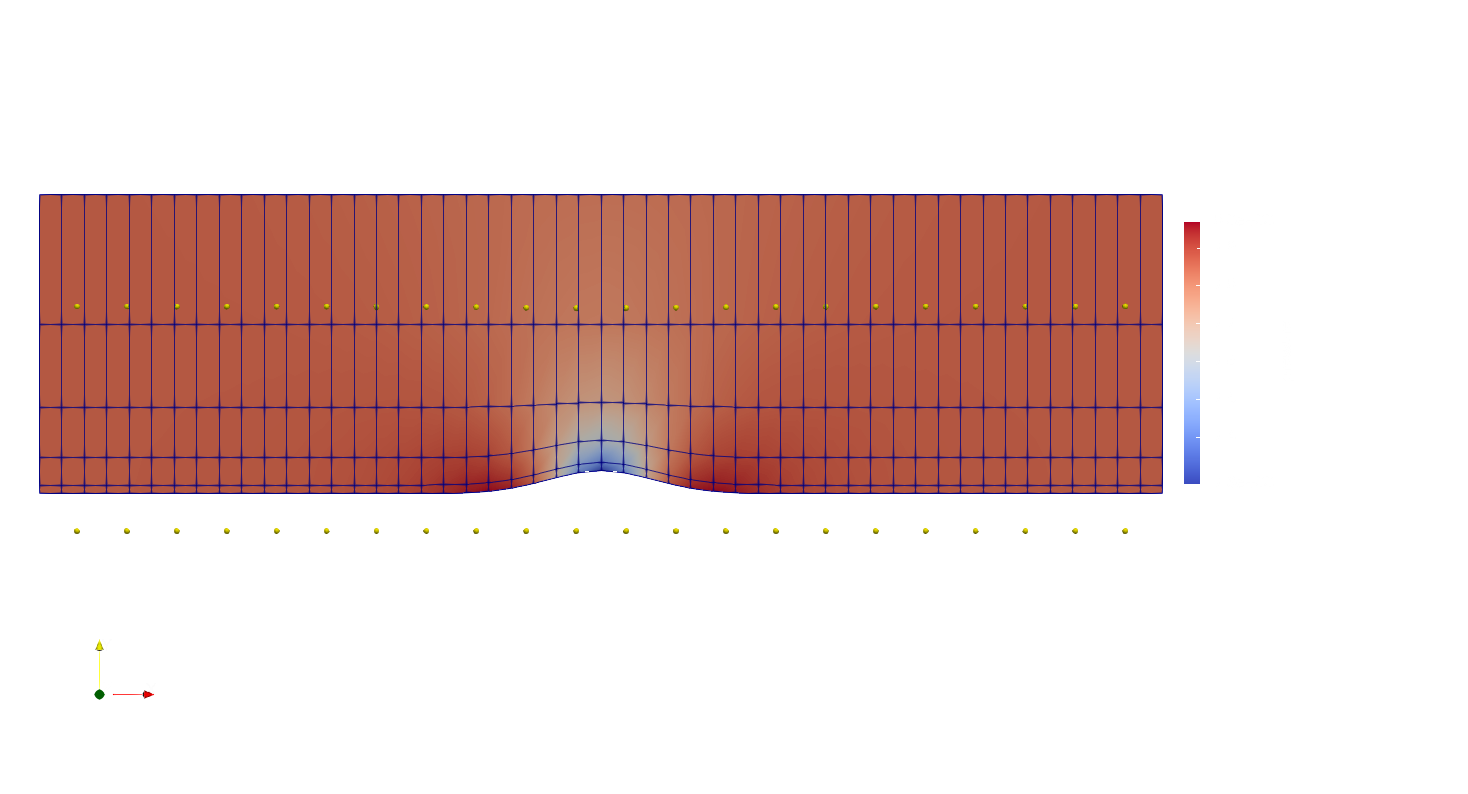}
        %}
    }
    \\
    \subfloat[][Target]{ \label{fig:target_gaussian_bump}
        \includegraphics
        [width=0.8\textwidth,trim={0.5cm 8.0cm 10.5cm 6cm},clip]
        {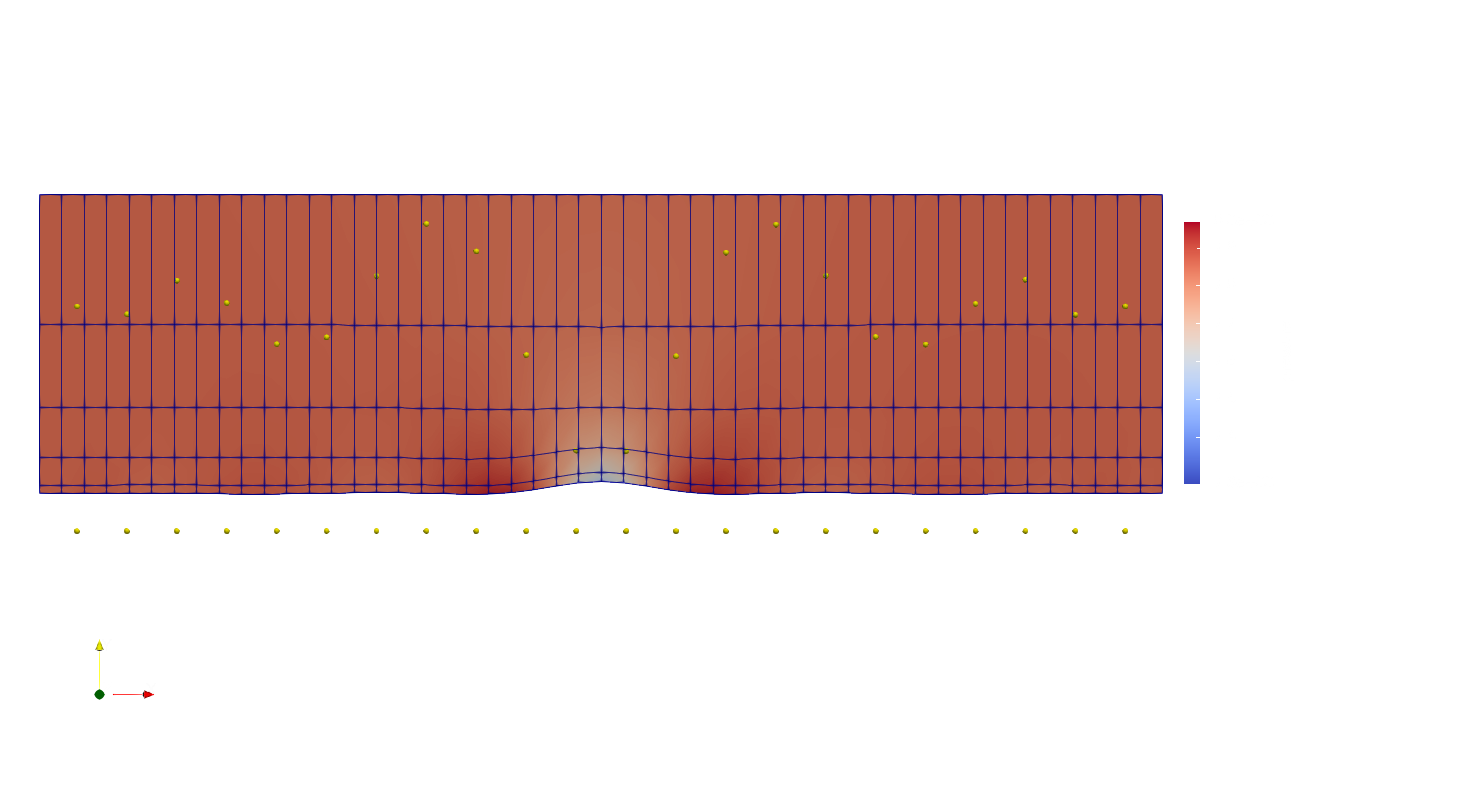}
    }
    \caption{Gaussian bump pressure distribution and FFD control points}
    \label{fig:gaussian_bump}
\end{figure}
The benchmark problem is the inverse target optimization of an inviscid channel flow from the High Order CFD workshop \cite{HOCFD2018}.
The channel is defined on $(x,y)=([-1.5,1.5],[0,0.8])$ and the lower surface follows a Gaussian bump
\begin{equation}
    y(x) = h \exp(-25x^2),
\end{equation}
where $h=0.0625$ for the initial grid and $h = 0.03125$ for the target grid as seen in Fig. \ref{fig:gaussian_bump}.
A FFD box is defined around the lower surface using $\mbf{Q} = (-1.4,-0.1)$, $\mbs{S} = (2.8, 0)$, and $\mbs{T} = (0,0.6)$.
The control variables are the $y$-component of the FFD control points, excluding the left-most, right-most, and bottom points to ensure a valid grid.

We impose total pressure and total temperature at the inlet, as well as an inflow Mach number of 0.3, while a static pressure ratio is prescribed at the subsonic outlet.
Characteristic boundary conditions as described in \cite{Carlson2011} are imposed implicitly.
The objective function is defined as
\begin{equation}
    \objective(\state,\des) = 
    \int_{\boundary} 
    \sum_i \left(
        u_i - u_i^{t}
    \right)^2 \d \boundary
    ,
\end{equation}

\subsection{Verification}

Since the test case is subsonic and the geometry is smooth the flow is isentropic.
The $L_2$-norm of the entropy error is used to measure the solution accuracy.
\begin{equation}
    E_{\text{entropy}} = \sqrt{\int_{\domain} \left( \frac{p}{p_\infty} \left(\frac{\rho_\infty}{\rho} \right)^\gamma - 1 \right)^2 \d \domain}.
\end{equation}
The entropy error should converge at a rate of $p+1$ for a DG discretization of order $p$ as seen in Table \ref{tab:bump_conv_rate}.
\input{table/bump_convergence_table}

\subsection{Optimization results}

The inverse design problem is solved using the full-space with various preconditioners, the reduced-space Newton, and the reduced-space quasi-Newton for various polynomial degrees ($p=1,2,3$, corresponding to $m=4000, 9000, 16000$ discretization degrees of freedom) and for various numbers of control variables ($n=20, 40, 60, 80, 100$).
%{\color{blue} Note that $m$ corresponds to the number of degrees of freedom of the discretization.}
The entire set of results is concatenated in Table \ref{tab:big_table}.
%The full-space with various preconditioners are compared in terms of iterations, subiterations, and work.
%The best performing full-space framework is then compared against the reduced-space Newton and quasi-Newton.

%\input{table/summary_table}
\input{table/fixed_summary_table}

%\subsubsection{Full-space}

\begin{figure}[h]
    \centering
    \includegraphics[width=\linewidth, page=23]
    {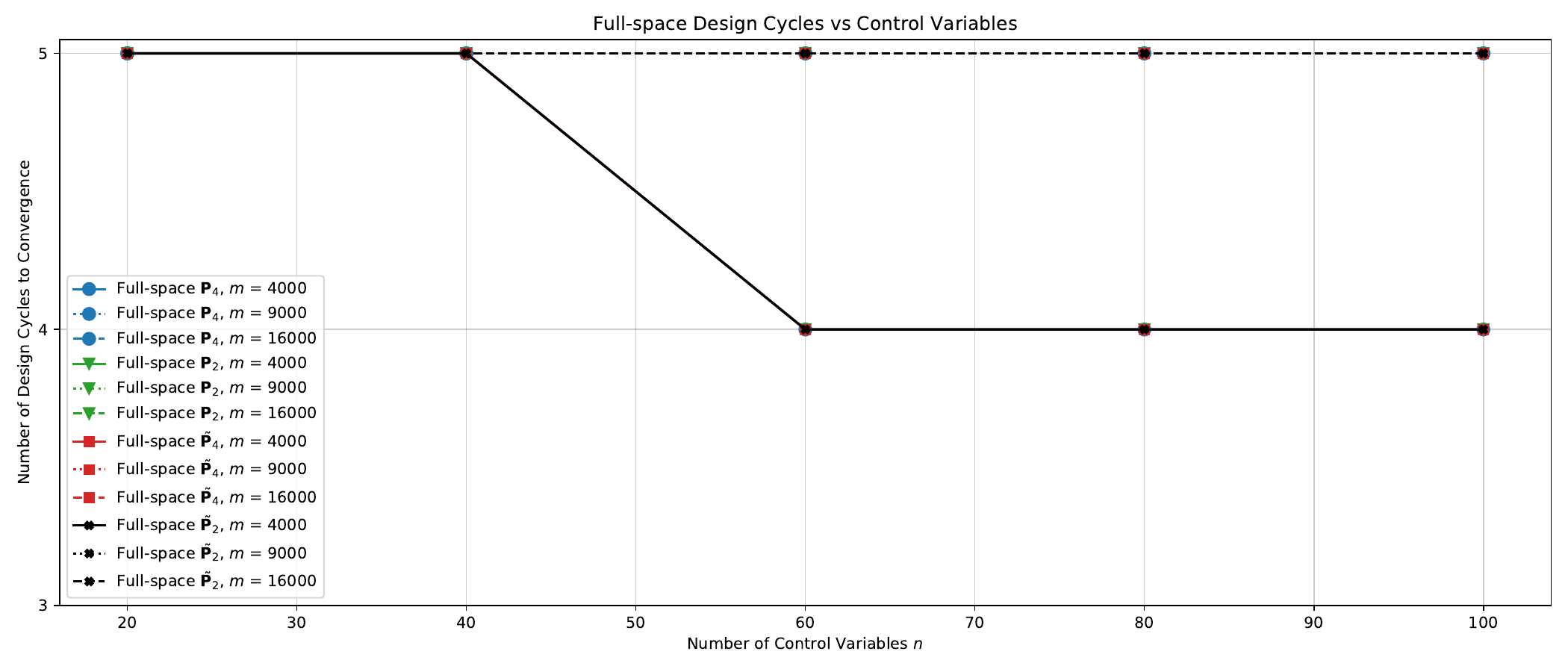}
    \caption{Full-space method convergence is independent of control variables size.}
    \label{fig:full_space_grad_it}
\end{figure}

The full-space framework is Newton-based and therefore converges independently of the control space size as seen in Fig. \ref{fig:full_space_grad_it}.
Since the same KKT system is solved to the same tolerance for various preconditioners, the search directions should be relatively similar.
Slight differences may occur since the preconditioner type will affect the accuracy of the various search direction components.
For example, it is expected that the $\mbf{P}_4$ preconditioner more accurately resolves the state and adjoint search directions, thus allowing a higher residual error in the control search direction.
Conversely, the $\Tilde{\mbf{P}}_2$ would distribute the allowable error amongst all three directions.

\begin{figure}[h]
    \centering
    \includegraphics[width=\linewidth, page=22]
    {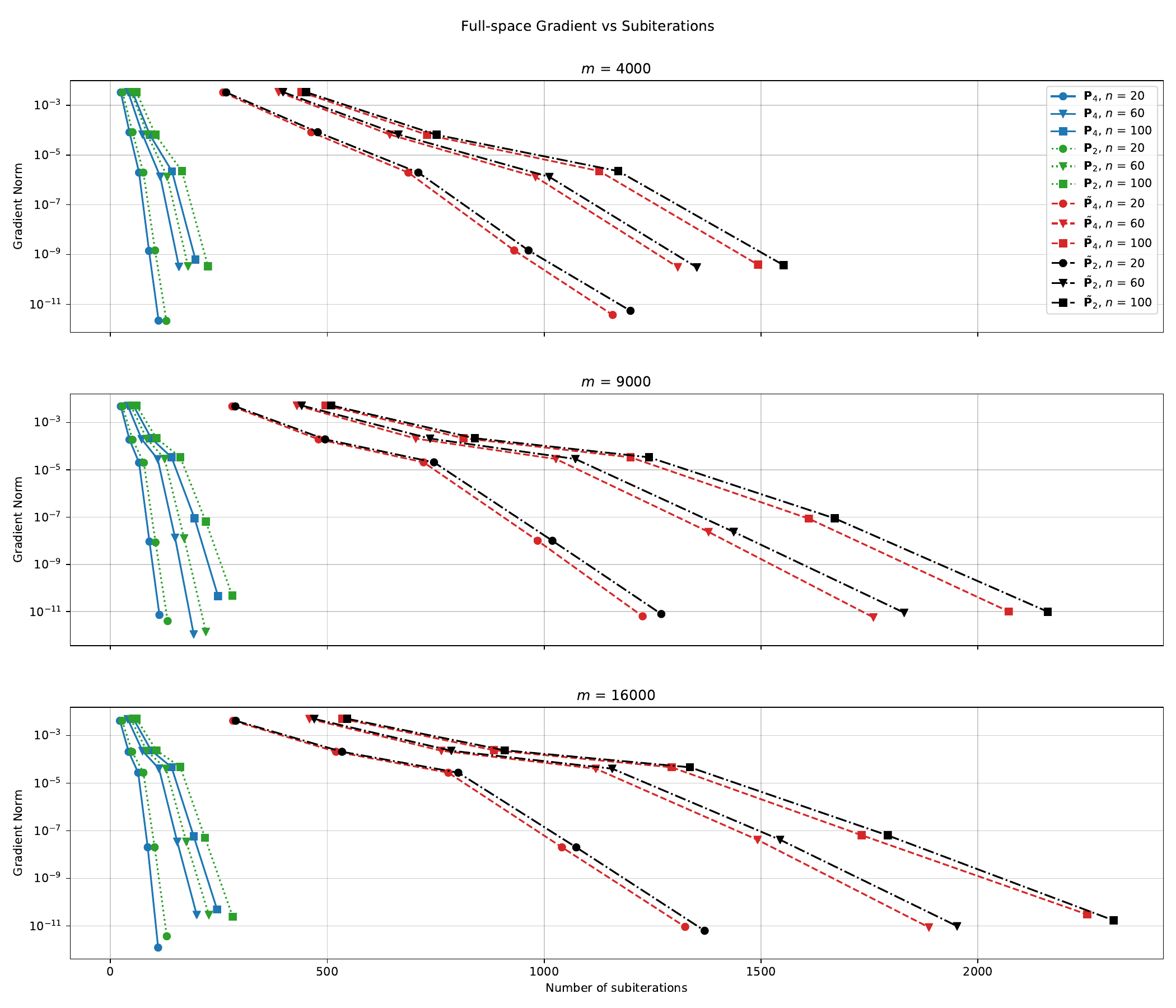}
    \caption{
    Effectiveness of preconditioners on the number of FGMRES iterations required to solve the KKT system.
    %Although we only plot the convergence for a spatial discretization of order 2, a similar trend occurs for P=1 and P=2.
    }
    \label{fig:full_space_grad_subit}
\end{figure}

\begin{figure}[h]
    \centering
    \includegraphics[width=\linewidth, page=25]
    {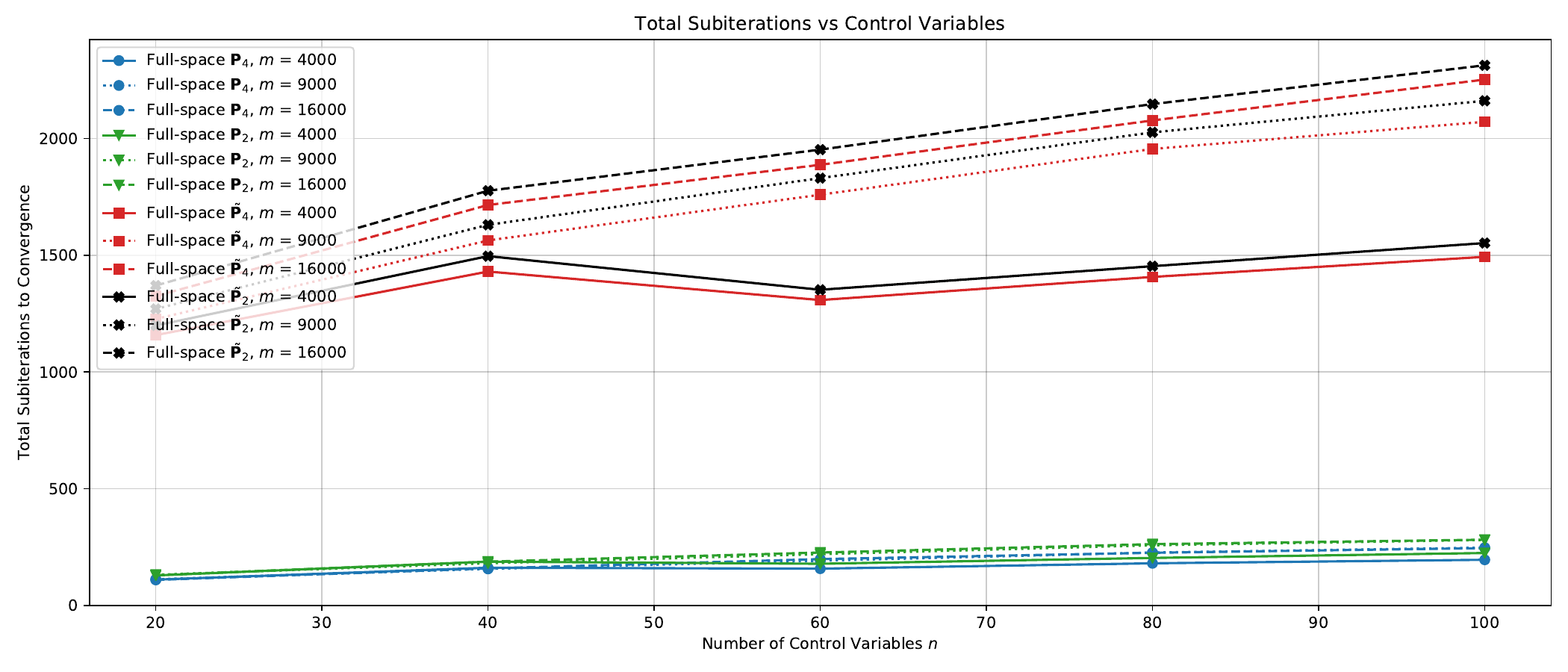}
    \caption{Full-space method total number of subiterations required to converge.}
    \label{fig:full_space_subit_ctl}
\end{figure}

The effectiveness of each preconditioner is observed in Fig. \ref{fig:full_space_grad_subit} \& \ref{fig:full_space_subit_ctl}.
As expected, the $\PFour$ and $\PTwo$ preconditioners greatly reduces the number of Krylov iterations since the eigenvalue spectrum is composed of at most $(\ndes+1)$ different values.
In fact, Table \ref{tab:big_table} shows that increasing the number of state variables does not affect the average number of subiterations, which confirms that the eigenvalue spectra of those preconditioned systems only depends on the control variables.
On the other hand, $\PFourA$ and $\PTwoA$ requires many more subiterations, but still allows FGMRES to converge, whereas the unpreconditioned system solve simply stalls.

\begin{figure}[h]
    \centering
    \includegraphics[width=\linewidth, page=24]
    {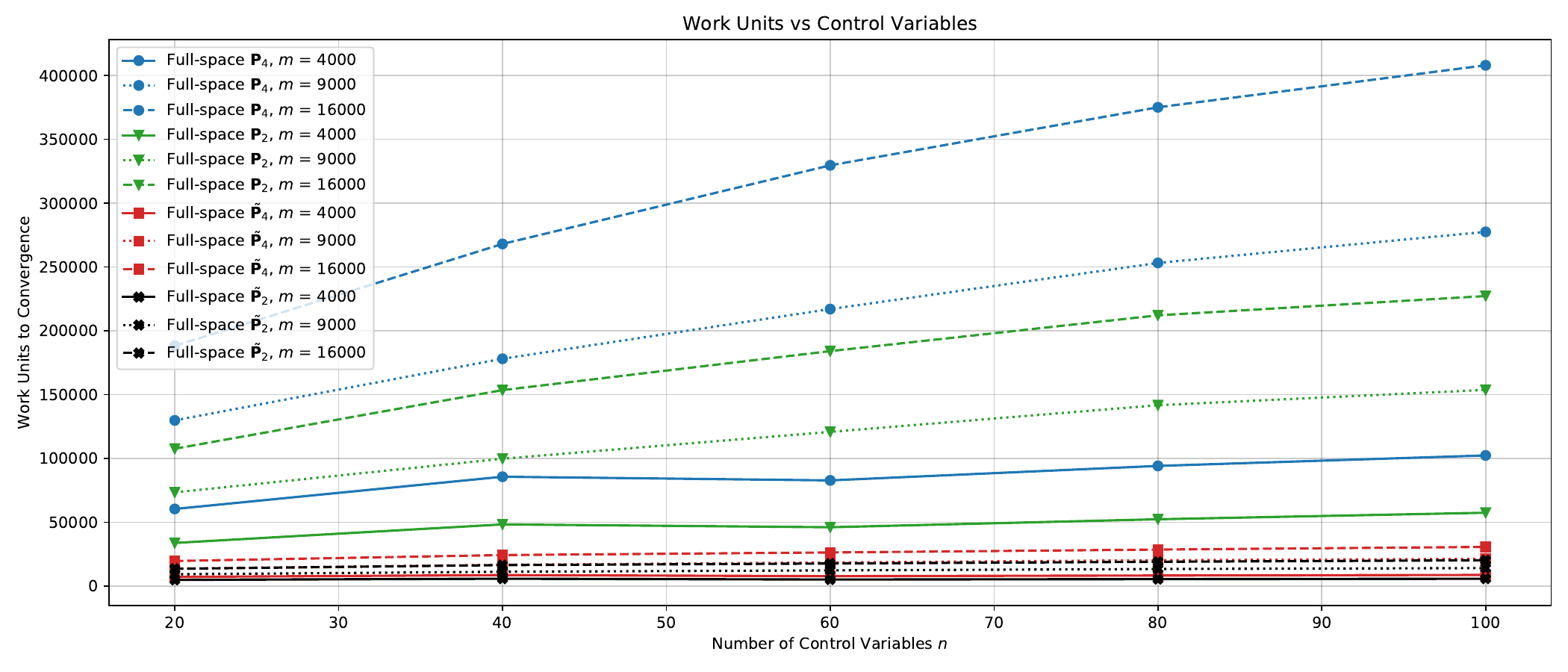}
    \caption{Full-space method work required to converge.}
    \label{fig:full_space_work_ctl}
\end{figure}

However, each application of $\PFour^{-1}$ and $\PTwo^{-1}$ is much costlier than $\PFourA^{-1}$ and $\PTwoA^{-1}$.
The total amount of work required to converge is shown in Fig. \ref{fig:full_space_work_ctl}.
As the control space or the state space increases in size, the cost of the $\PFour$ preconditioned full-space, that excelled in subiterations, shows a significant increase in cost.
%The opposite happens when $\PTwoA$ is used, where the work slowly increases as $\nstate$ or $\ndes$ increases although it required the most subiterations.
Table \ref{tab:big_table} shows that when $\PTwoA$ is used, the work slowly increases as $\nstate$ or $\ndes$ increases although it required the most subiterations.

\begin{figure}[h]
    \centering
    \includegraphics[width=\linewidth, page=27]
    {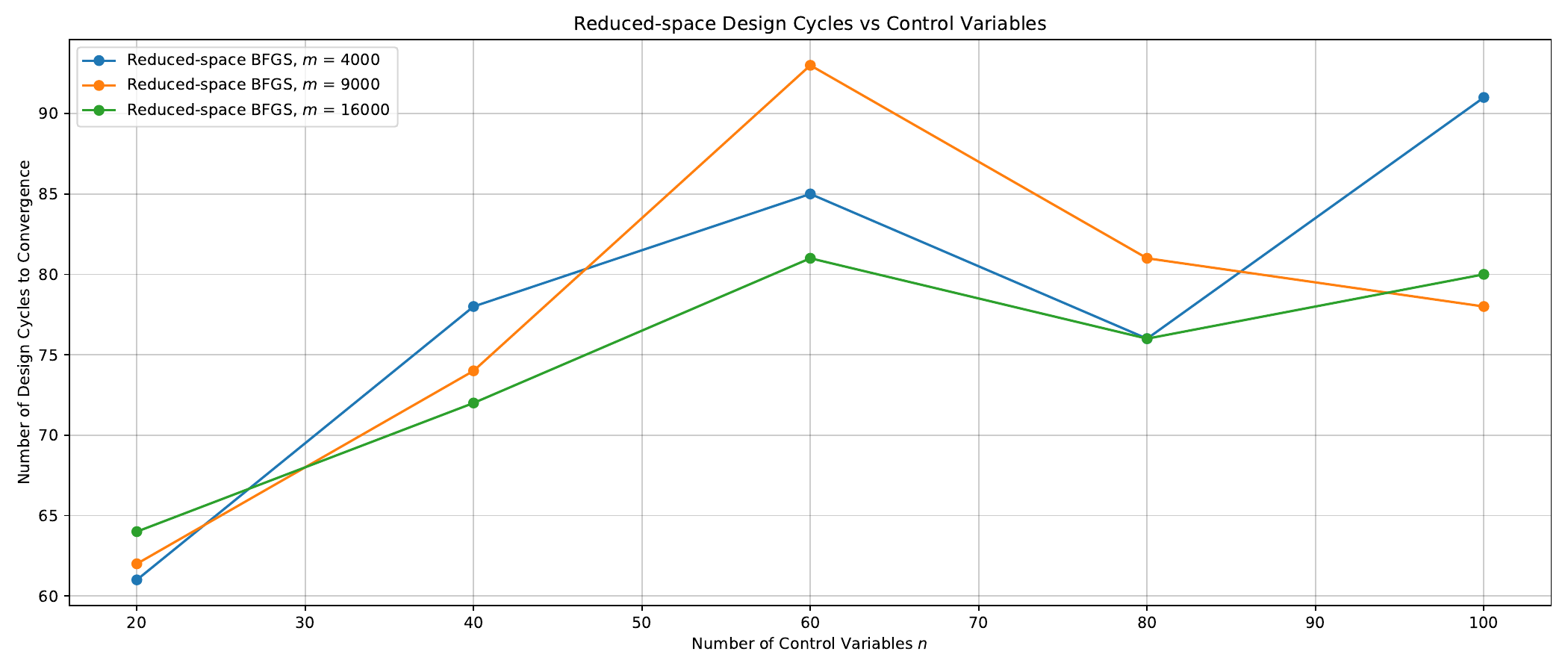}
    \caption{Reduced-space BFGS design cycles requires to converge. Note that the reduced-space Newton method requires a constant 4 design cycles regardless of $\ndes$ or $\nstate$.}
    \label{fig:bfgs_cycles_nctl}
\end{figure}

%\subsubsection{Reduced-space}

\begin{figure}[h]
    \centering
    \includegraphics[width=\linewidth, page=28]
    {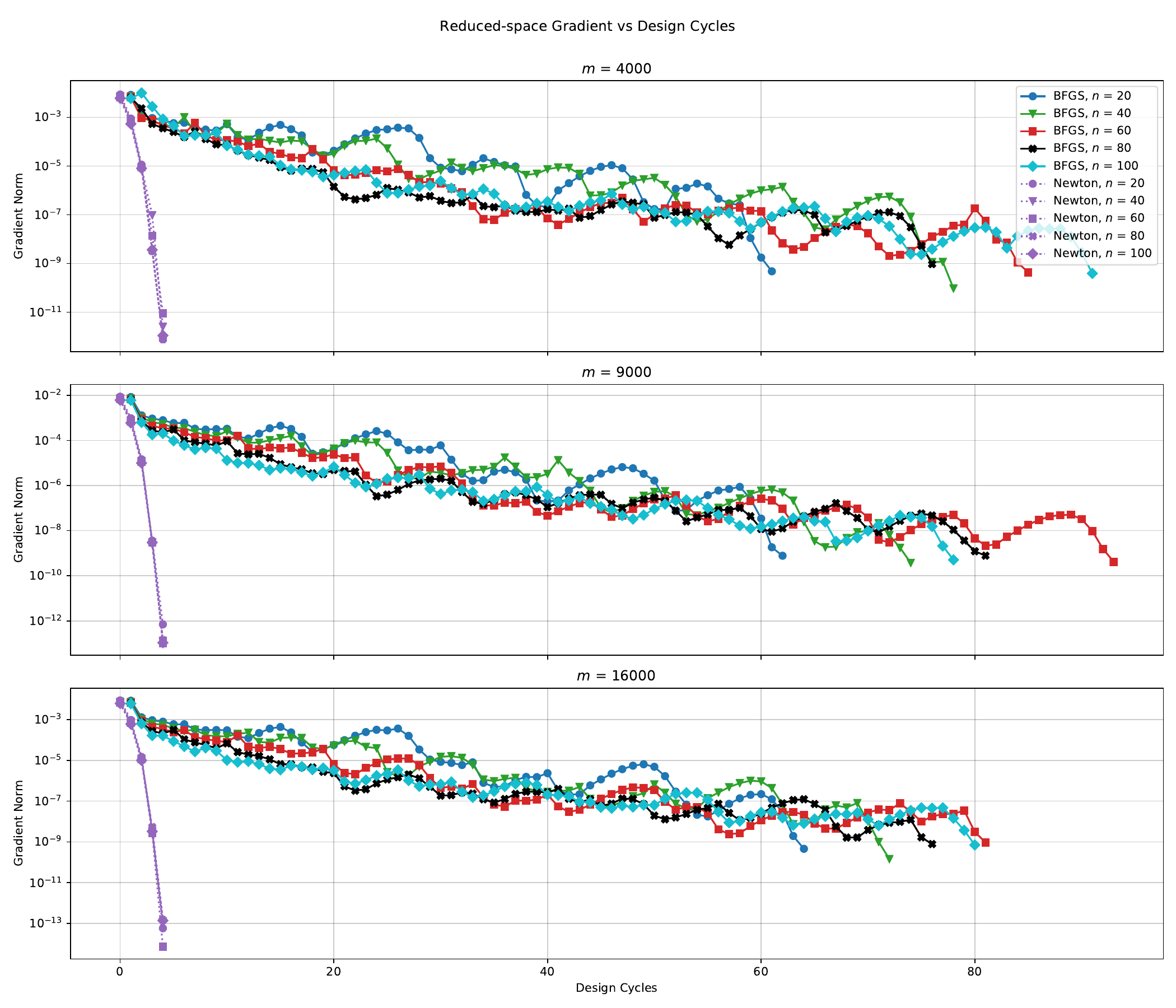}
    \caption{Reduced-space gradient convergence with design cycles.}
    \label{fig:reduced_grad_it}
\end{figure}

\begin{figure}[h]
    \centering
    \includegraphics[width=\linewidth, page=26]
    {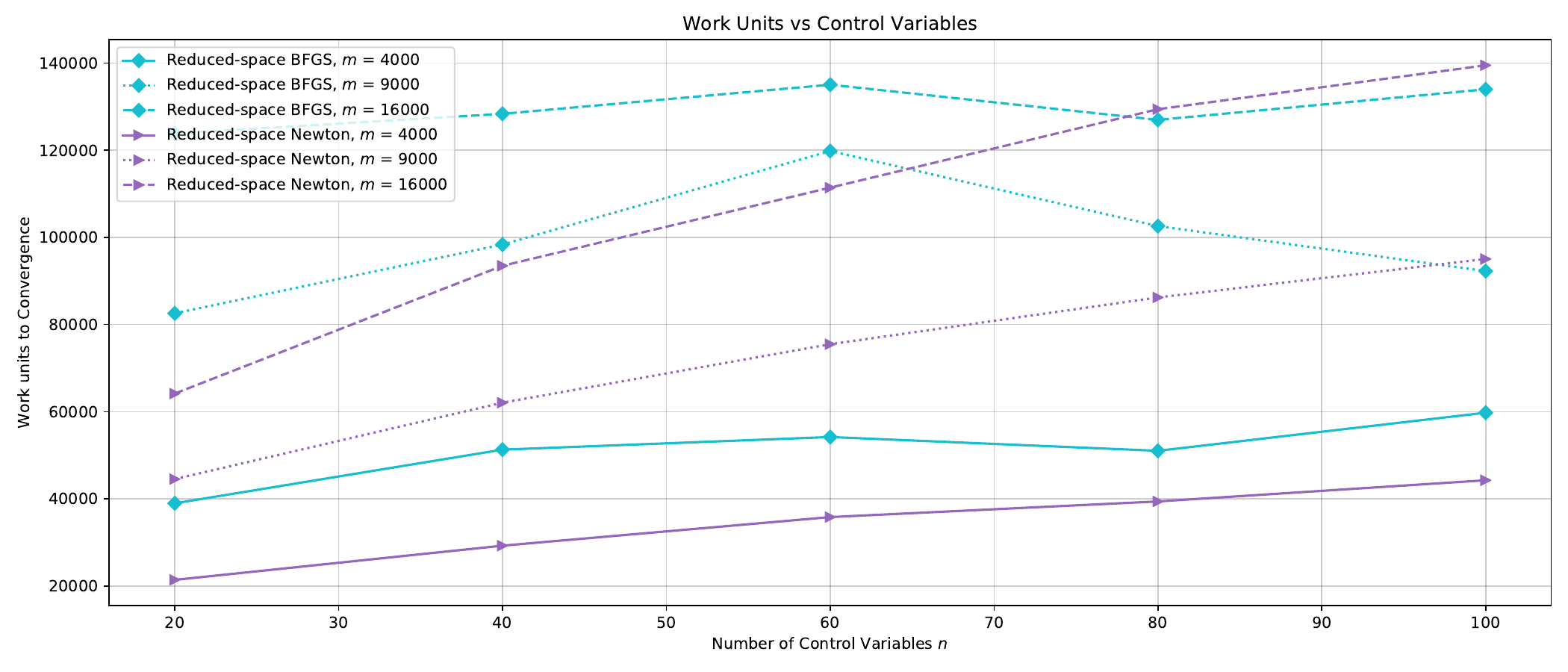}
    \caption{Work comparison between reduced-space approaches.}
    \label{fig:reduced_work_ctl}
\end{figure}

%The quasi-Newton compared against the Newton method follows a similar pattern between comparing the $\PFour$ and $\PTwo$ preconditioners and their $\PFourA$ and $\PTwoA$ counterpart as seen in Fig. \ref{fig:reduced_grad_it}.
The quasi-Newton is compared against the Newton method in Fig. \ref{fig:reduced_grad_it}.
The Newton method converges independently of the control space size just as the full-space method, whereas the BFGS method requires an increasing number of design cycles as seen Fig. \ref{fig:bfgs_cycles_nctl}.
However, Table \ref{tab:big_table} shows that the number of subiterations within the Newton's method linear solve increases as the eigenvalue spectrum widens.
Looking at the total work to convergence in Fig. \ref{fig:reduced_work_ctl} gives a clearer picture.
%\st{The quasi-Newton method requires significantly more work units than the Newton method because of multiple flow solves required.}
The reduced-space Newton method outperforms the quasi-Newton approach for lower numbers of control variables, but its cost becomes comparable around 100 control variables.
%Both reduced-space methods have comparable costs at 100 control variables,

\begin{figure}[h]
    \centering
    \includegraphics[width=\linewidth, page=29]
    {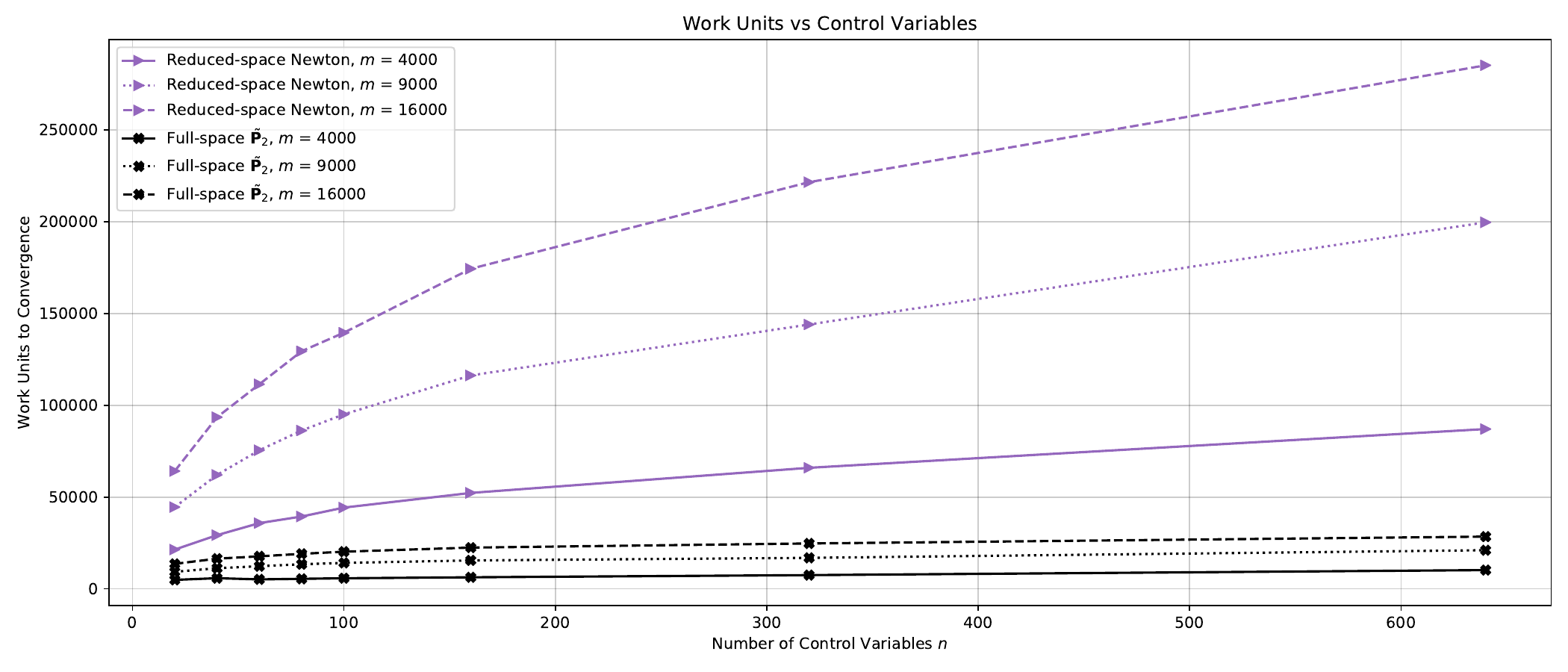}
    \caption{Scaling of work required to converge with control variables between reduced and full-space approaches.}
    \label{fig:fullreduced_work_ctl}
\end{figure}

Finally, we compare the reduced-space approach with the best performing full-space approach in Fig. \ref{fig:fullreduced_work_ctl}.
Additional data points at $n = 160, 320, 640$ have been evaluated to observe the scaling of both algorithms.
Unfortunately, the reduced-space BFGS algorithm fails to find a suitable direction during a linesearch for larger number of control variables.
% \& \ref{fig:fullreduced_work_sim}.
The $\PTwoA$ preconditioned full-space method outperforms the reduced-space Newton method for every design set.
Furthermore, by inspecting the average work per cycle of Table \ref{tab:big_table}, we also see that the full-space method scales better than the reduced-space Newton when the number of state variables $n$ increases.

However, the increase in the full-space average subiterations per cycle shows that a simple BFGS preconditioner is not sufficient to achieve control-independent convergence.
Since the full-space only updates its reduced-Hessian approximation 4-5 times, it is expected that the identity-initialized BFGS preconditioner becomes unsuitable for a large set of control variables.
On the other hand, the reduced-space quasi-Newton is able to benefit from more BFGS updates by performing a larger number of less expensive design cycles.

%\subsubsection{Full-space versus reduced-space}

\begin{figure}[h]
    \centering
    \includegraphics[width=\linewidth, page=30]
    {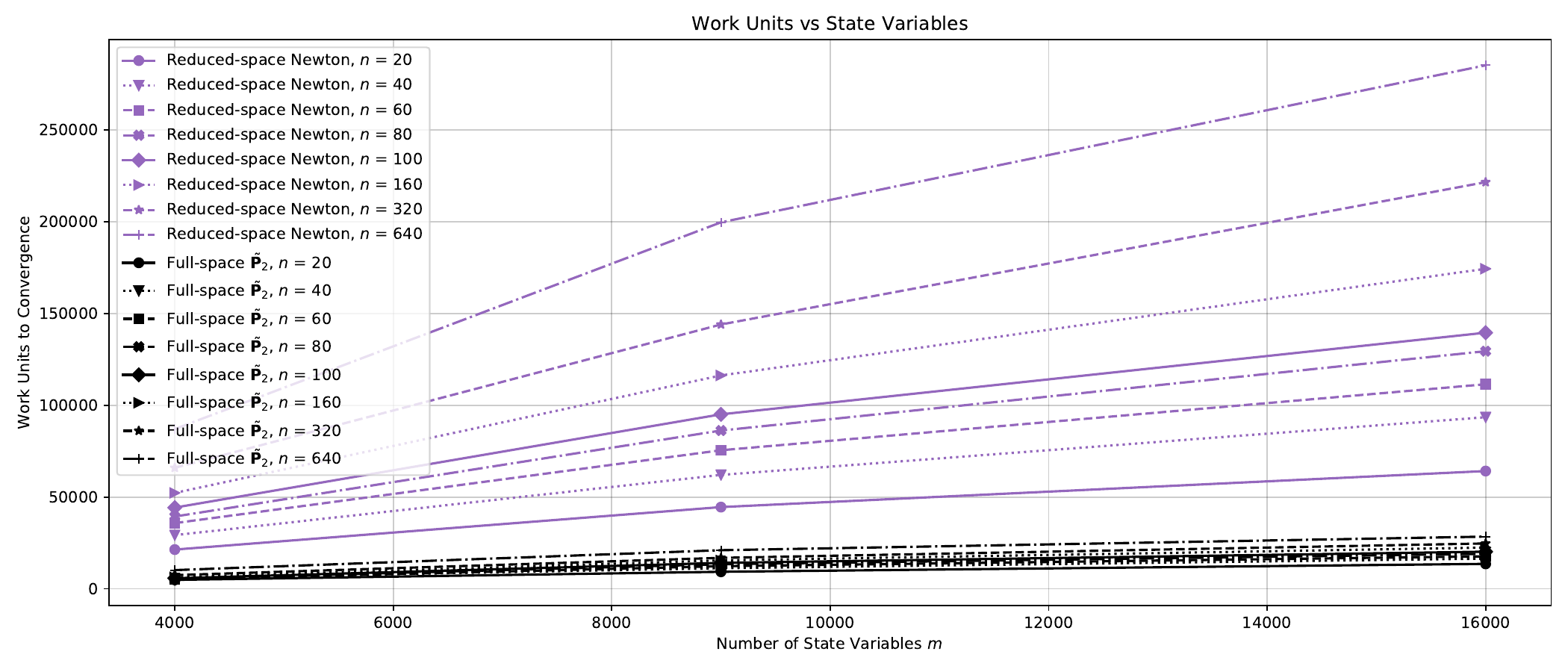}
    \caption{Scaling of work required to converge with state variables between reduced and full-space approaches.}
    \label{fig:fullreduced_work_sim}
\end{figure}

%% file: table/bump_convergence_table.tex
\begin{table}[h]
\centering
\resizebox{0.65\textwidth}{!}
{%
\begin{tabular}{r S[table-format=7.0] S[table-format=7.0] ccr}
\toprule
\textbf{Degree $p$} 
& \textbf{Cells }
& \textbf{Number of DoFs }
& \textbf{$L_2$ entropy error }
& \textbf{Rate} \\
\midrule
\midrule
1                    & 64                   & 1024                        & 5.73\Eminus{03}             & -                    \\
1                    & 256                  & 4096                        & 1.30\Eminus{03}             & 2.14                 \\
1                    & 1024                 & 16384                       & 2.26\Eminus{04}             & 2.53                 \\
1                    & 4096                 & 65536                       & 4.30\Eminus{05}             & 2.39                 \\
1                    & 16384                & 262144                      & 9.07\Eminus{06}             & 2.25                 \\
\midrule
2                    & 64                   & 2304                        & 9.86\Eminus{04}             & -                    \\
2                    & 256                  & 9216                        & 7.68\Eminus{05}             & 3.68                 \\
2                    & 1024                 & 36864                       & 7.19\Eminus{06}             & 3.42                 \\
2                    & 4096                 & 147456                      & 7.91\Eminus{07}             & 3.18                 \\
2                    & 16384                & 589824                      & 9.34\Eminus{08}             & 3.08                 \\
\midrule
3                    & 64                   & 4096                        & 1.57\Eminus{04}             & -                    \\
3                    & 256                  & 16384                       & 1.03\Eminus{05}             & 3.94                 \\
3                    & 1024                 & 65536                       & 6.80\Eminus{07}             & 3.91                 \\
3                    & 4096                 & 262144                      & 5.07\Eminus{08}             & 3.75                 \\
3                    & 16384                & 1048576                     & 3.23\Eminus{09}             & 3.97                 \\
\bottomrule
\end{tabular}
}
\caption{Orders of accuracy table demonstrating $L_2$ norm of the entropy error diminishing at a rate of of $p+1$ for a discretization of order $p$.}
\label{tab:bump_conv_rate}
\end{table}

%% file: table/fixed_summary_table.tex
\begin{table}[ht]
    \centering
\resizebox{\textwidth}{!}
{%
\begin{tabular}{
r % $\ndes$ 
%c % Cycles 
S[table-format=2.0] @{\,\,\( / \)\,}  S[table-format=3.0] @{\,\,\( / \)\,} S[table-format=3.0]
%c % Total Subits. 
S[table-format=4.0] @{\,\,\( / \)\,}  S[table-format=4.0] @{\,\,\( / \)\,} S[table-format=4.0]
%c % Subits./cycle 
S[table-format=4.1] @{\,\,\( / \)\,}  S[table-format=4.1] @{\,\,\( / \)\,} S[table-format=4.1]
%c % Total Work 
S[table-format=6.0] @{\,\,\( / \)\,}  S[table-format=6.0] @{\,\,\( / \)\,} S[table-format=6.0]
%c % Work/cycle
S[table-format=6.0] @{\,\,\( / \)\,}  S[table-format=6.0] @{\,\,\( / \)\,} S[table-format=6.0]
}
\toprule 
\textbf{$\ndes$ }
& \multicolumn{3}{c}{\textbf{Cycles }}
& \multicolumn{3}{c}{\textbf{Total Subits. }}
& \multicolumn{3}{c}{\textbf{Avg. Subits./Cycle }}
& \multicolumn{3}{c}{\textbf{Total Work }}
& \multicolumn{3}{c}{\textbf{Avg. Work/Cycle} } \\
\midrule
\midrule 
\multicolumn{16}{c}{\textbf{Full-space $\PFour$}} \\ 
\midrule 
   20  &   5 &    5 &    5 & 111 &  113 &  110 & 22.2 &  22.6 &  22.0 &  60453 &  129849 &  188277 &  12091 &   25970 &   37655 \\ 
   40  &   5 &    5 &    5 & 162 &  157 &  159 & 32.4 &  31.4 &  31.8 &  85696 &  178114 &  268069 &  17139 &   35623 &   53614 \\ 
   60  &   4 &    5 &    5 & 158 &  192 &  199 & 39.5 &  38.4 &  39.8 &  82837 &  217048 &  329648 &  20709 &   43410 &   65930 \\ 
   80  &   4 &    5 &    5 & 181 &  227 &  226 & 45.2 &  45.4 &  45.2 &  94181 &  253188 &  375115 &  23545 &   50638 &   75023 \\ 
  100  &   4 &    5 &    5 & 196 &  248 &  246 & 49.0 &  49.6 &  49.2 & 102381 &  277537 &  408023 &  25595 &   55507 &   81605 \\ 
\midrule 
\multicolumn{16}{c}{\textbf{Full-space $\PTwo$}} \\ 
\midrule 
   20  &   5 &    5 &    5 & 129 &  132 &  130 & 25.8 &  26.4 &  26.0 &  33828 &   73457 &  107643 &   6766 &   14691 &   21529 \\ 
   40  &   5 &    5 &    5 & 188 &  181 &  188 & 37.6 &  36.2 &  37.6 &  48347 &   99804 &  153567 &   9669 &   19961 &   30713 \\ 
   60  &   4 &    5 &    5 & 179 &  220 &  227 & 44.8 &  44.0 &  45.4 &  46129 &  120854 &  184120 &  11532 &   24171 &   36824 \\ 
   80  &   4 &    5 &    5 & 204 &  258 &  263 & 51.0 &  51.6 &  52.6 &  52391 &  141722 &  212156 &  13098 &   28344 &   42431 \\ 
  100  &   4 &    5 &    5 & 225 &  281 &  282 & 56.2 &  56.2 &  56.4 &  57491 &  153734 &  227260 &  14373 &   30747 &   45452 \\ 
\midrule 
\multicolumn{16}{c}{\textbf{Full-space $\PFourA$}} \\ 
\midrule 
   20  &   5 &    5 &    5 &1158 & 1227 & 1325 &231.6 & 245.4 & 265.0 &   7223 &   13634 &   19695 &   1445 &    2727 &    3939 \\ 
   40  &   5 &    5 &    5 &1430 & 1563 & 1715 &286.0 & 312.6 & 343.0 &   8640 &   16710 &   24329 &   1728 &    3342 &    4866 \\ 
   60  &   4 &    5 &    5 &1308 & 1759 & 1887 &327.0 & 351.8 & 377.4 &   7800 &   18505 &   26372 &   1950 &    3701 &    5274 \\ 
   80  &   4 &    5 &    5 &1407 & 1955 & 2077 &351.8 & 391.0 & 415.4 &   8316 &   20299 &   28627 &   2079 &    4060 &    5725 \\ 
  100  &   4 &    5 &    5 &1493 & 2071 & 2252 &373.2 & 414.2 & 450.4 &   8764 &   21361 &   30709 &   2191 &    4272 &    6142 \\ 
\midrule 
\multicolumn{16}{c}{\textbf{Full-space $\PTwoA$}} \\ 
\midrule 
   20  &   5 &    5 &    5 &1199 & 1270 & 1370 &239.8 & 254.0 & 274.0 &   4893 &    9296 &   13605 &    979 &    1859 &    2721 \\ 
   40  &   5 &    5 &    5 &1496 & 1630 & 1776 &299.2 & 326.0 & 355.2 &   5811 &   11250 &   16466 &   1162 &    2250 &    3293 \\ 
   60  &   4 &    5 &    5 &1352 & 1830 & 1952 &338.0 & 366.0 & 390.4 &   5161 &   12336 &   17706 &   1290 &    2467 &    3541 \\ 
   80  &   4 &    5 &    5 &1453 & 2026 & 2147 &363.2 & 405.2 & 429.4 &   5473 &   13400 &   19077 &   1368 &    2680 &    3815 \\ 
  100  &   4 &    5 &    5 &1552 & 2161 & 2313 &388.0 & 432.2 & 462.6 &   5779 &   14133 &   20249 &   1445 &    2827 &    4050 \\ 
  160  &   4 &    5 &    5 &1721 & 2409 & 2628 &430.2 & 481.8 & 525.6 &   6302 &   15480 &   22469 &   1575 &    3096 &    4494 \\ 
  320  &   4 &    5 &    5 &2104 & 2668 & 2948 &526.0 & 533.6 & 589.6 &   7486 &   16886 &   24723 &   1871 &    3377 &    4945 \\ 
  640  &   5 &    5 &    5 &2920 & 3435 & 3477 &584.0 & 687.0 & 695.4 &  10212 &   21050 &   28451 &   2042 &    4210 &    5690 \\ 
\midrule 
\multicolumn{16}{c}{\textbf{Reduced-space Newton}} \\ 
\midrule 
   20  &   4 &    4 &    4 &  68 &   67 &   65 & 17.0 &  16.8 &  16.2 &  21397 &   44526 &   64156 &   5349 &   11131 &   16039 \\ 
   40  &   4 &    4 &    4 & 100 &  101 &  103 & 25.0 &  25.2 &  25.8 &  29248 &   62060 &   93477 &   7312 &   15515 &   23369 \\ 
   60  &   4 &    4 &    4 & 128 &  126 &  126 & 32.0 &  31.5 &  31.5 &  35809 &   75482 &  111420 &   8952 &   18871 &   27855 \\ 
   80  &   4 &    4 &    4 & 142 &  147 &  149 & 35.5 &  36.8 &  37.2 &  39397 &   86220 &  129417 &   9849 &   21555 &   32354 \\ 
  100  &   4 &    4 &    4 & 162 &  163 &  161 & 40.5 &  40.8 &  40.2 &  44256 &   95040 &  139496 &  11064 &   23760 &   34874 \\ 
  160  &   4 &    4 &    4 & 193 &  202 &  206 & 48.2 &  50.5 &  51.5 &  52244 &  116218 &  174302 &  13061 &   29055 &   43575 \\ 
  320  &   4 &    4 &    4 & 246 &  253 &  264 & 61.5 &  63.2 &  66.0 &  65920 &  143962 &  221465 &  16480 &   35991 &   55366 \\ 
  640  &   4 &    4 &    4 & 336 &  359 &  346 & 84.0 &  89.8 &  86.5 &  87027 &  199590 &  285098 &  21757 &   49898 &   71274 \\ 
\midrule 
  \multicolumn{16}{c}{\textbf{Reduced-space BFGS}} \\ 
\midrule 
   20  &  61 &   62 &   64 &   0 &    0 &    0 &  0.0 &   0.0 &   0.0 &  38967 &   82569 &  123860 &    639 &    1332 &    1935 \\ 
   40  &  78 &   74 &   72 &   0 &    0 &    0 &  0.0 &   0.0 &   0.0 &  51288 &   98358 &  128359 &    658 &    1329 &    1783 \\ 
   60  &  85 &   93 &   81 &   0 &    0 &    0 &  0.0 &   0.0 &   0.0 &  54183 &  119824 &  135045 &    637 &    1288 &    1667 \\ 
   80  &  76 &   81 &   76 &   0 &    0 &    0 &  0.0 &   0.0 &   0.0 &  51022 &  102596 &  126967 &    671 &    1267 &    1671 \\ 
  100  &  91 &   78 &   80 &   0 &    0 &    0 &  0.0 &   0.0 &   0.0 &  59764 &   92270 &  133987 &    657 &    1183 &    1675 \\ 
\bottomrule
\end{tabular}
}
    \caption{Convergence summary of the various methods discussed. $\left(\cdot\right) / \left(\cdot\right) / \left(\cdot\right)$ are the data points given for $m = 4000, 9000, 16000$. }
\label{tab:big_table}
\end{table}

%% file: 9_conclusionandacknowledgements.tex
\section{Conclusions}

The necessary ingredients to implement the full-space LNKS approach within an ASO framework have been laid out, and a cost analysis of the various derivative operators is presented. A re-derivation of the $\mbf{P}_2$ and $\mbf{P}_4$ preconditioners of~\cite{Biros2005a} yielded new preconditioners that differ from the original.
%\sout{The preconditioned systems structure was shown to be different than 
%previously stated from~\cite{Biros2005a}.}
However, it is analytically and numerically shown that the eigenvalue spectra are the same for the corrected preconditioned systems.
%Furthermore, the numerical results confirm that the $\PFour$ and $\PTwo$ preconditioned systems on

The full-and reduced-space approaches have been applied and compared on a benchmark aerodynamic problem.
A dimensional analysis shows that the full-space approach with $\PTwoA$ preconditioning is the most effective method.
However, its scaling will largely depend on the development of reduced-Hessian preconditioners and flow Jacobian preconditioners.
Future work should aim to achieve control-independent convergence through the development of reduced-Hessian preconditioners, such as approximating the initial BFGS Hessian \cite{ShiDong2018}.

While the cost analysis in the present work reflects our implementation, further savings through sum factorization, collocation, and other code optimizations will affect the work ratios and is subject to further investigation.
Additional work reduction within each design cycle can be achieved through inexact solves as proposed by~\cite{Biros2005} for the full-space approach,~\cite{Hicken2014} for the reduced-space Newton approach, and~\cite{Brown2017} for the reduced-space quasi-Newton approach.

%Finally, the adoption of full-space methods within ASO will depend on its ability to handle complex flow problems with more nonlinearities.
%As theTransitioning to full-space methods will require

\section*{Acknowledgements}

The authors gratefully acknowledge the generous support from 
the Natural Sciences and Engineering Research Council (NSERC) and the McGill Engineering Doctoral Award. We also thank Compute Canada for providing the computational facilities.

%% file: 9z_appendix.tex
\clearpage
\appendix

\onecolumn

\newcommand{\Rx}{{\residuals}_{\des}}
\newcommand{\Rs}{{\residuals}_{\state}}
\newcommand{\Lss}{\lagrangian_{\state \state}}
\newcommand{\Lsx}{\lagrangian_{\state \des}}
\newcommand{\Lxs}{\lagrangian_{\des \state}}
\newcommand{\Lxx}{\lagrangian_{\des \des}}

\newcommand{\Lyy}{\lagrangian_{yy}}
\newcommand{\Lzz}{\lagrangian_{zz}}

\newcommand{\KoneInv}{\left(\KKTmat_1^{-1}\right)}
\newcommand{\KtwoInv}{\left(\KKTmat_2^{-1}\right)}
\newcommand{\PfourInv}{\left(\mbf{P}_4^{-1}\right)}

\input{9z_precon_addendum}

\input{9z_precond_eig}

\input{9z_precond_proof}

\input{9z_flops}

%% file: 9z_precon_addendum.tex
\section{Preconditioned system corrections}
\label{app:preconditioner_inverses}
%\subsection{Abbreviations}
The purpose of this appendix is to correct the derivation from Biros and Ghattas' preconditioned system \cite{Biros2005a}.
Notation from Sec. \ref{sec:optimization_formulation} and \ref{sec:full_space_precond} is used in the following appendix.
Additionally, we will be using the following notation to simplify and shorten the derivation
\begin{align}
    \Lyy &=  \Lsx - \Lss \Rs^{-1} \Rx,
\\
    \Lyy^T &=  \Lxs -  \Rx^T \Rs^{-T} \Lss.
\end{align}
\begin{align}
    \Lzz &= \Rx^{T} \Rs^{-T} \Lss \Rs^{-1} \Rx - \Rx^{T} \Rs^{-T} \Lsx - \Lxs \Rs^{-1} \Rx + \Lxx
\\ & = 
    -\Rx^T \Rs^{-T} \Lyy
    - \Lxs \Rs^{-1} \Rx
    + \Lxx
\\ & = 
    - \Lyy^T \Rs^{-1} \Rx
    -\Rx^T \Rs^{-T} \Lsx
    + \Lxx
.
\end{align}

%\subsection{$\mathbf{P}_4$ preconditioning}
When $\mbf{P}_4^{-1}$ is applied to the KKT matrix $\KKTmat$, we obtain
\begin{equation}
\begin{split}
\mbf{P}_4^{-1} \KKTmat
%&=
%\begin{bmatrix}
%\PfourInv_{11} & \PfourInv_{12} & \PfourInv_{13} \\
%\PfourInv_{21} & \PfourInv_{22} & \PfourInv_{23} \\
%\PfourInv_{31} & \PfourInv_{32} & \PfourInv_{33}
%\end{bmatrix}
%\begin{bmatrix}
% \Lss & \Lsx & \Rs^T \\
% \Lxs & \Lxx & \Rx^T \\
% \Rs & \Rx & \bzero
%\end{bmatrix}
=
\begin{bmatrix}
\identity_{\state} 
&
\Rs^{-1}  \Rx
    \left(
    \identity_{\des} - 
    \approxHessian^{-1} \Lzz
    \right)
& \bzero \\
\bzero & \approxHessian^{-1} \Lzz & \bzero \\
\bzero 
&
  \Rs^{-T} \Lyy
  \left(
  \identity_{\des} - \approxHessian^{-1} \Lzz
  \right)
&
\identity_{\state}
\end{bmatrix}
\label{eq:p4invK}
\end{split}
,
\end{equation}
which differs from the original derivation in \cite{Biros2005a} that states
\begin{equation}
\begin{split}
\left(\mbf{P}_4^{-1} \KKTmat \right)_{\text{BG}}
=
\begin{bmatrix}
\identity_{\state} 
&
\bzero
& \bzero \\
\bzero & \approxHessian^{-1} \Lzz & \bzero \\
\bzero  & \bzero & \identity_{\state}
\end{bmatrix}
\end{split}
.
\label{eq:p4invK_error}
\end{equation}

Furthermore, when $\mbf{P}_2^{-1}$ applied to the KKT matrix $\KKTmat$, we obtain
\begin{equation}
\begin{split}
\mbf{P}_2^{-1} \KKTmat
=
\begin{bmatrix}
 \identity_{\state}
 -
 \Rs^{-1} \Rx
 \approxHessian^{-1} \Lyy^T
&
% \Rs^{-1} \Rx
% \left(
%     \identity_{\des}
%     - \approxHessian^{-1}  \Lzz
% \right)
% -
% \Rs^{-1} \Rx
% \approxHessian^{-1}  \Lyy^T
% \Rs^{-1} \Rx
% 
 \Rs^{-1} \Rx
 - \Rs^{-1} \Rx \approxHessian^{-1}  
 \left(
 \Lzz 
 -
 \Lyy^T \Rs^{-1} \Rx
 \right)
&
\bzero
\\
 \approxHessian^{-1} \Lyy^T
&
 \approxHessian^{-1} 
 \left(
 \Lzz 
 + \Lyy^T \Rs^{-1} \Rx
 \right)
&
\bzero
\\
 \Rs^{-T}  \Lss
&
 \Rs^{-T} 
 \Lsx
&
\identity_{\state}
\end{bmatrix}
\end{split}
,
\label{eq:p2invK}
\end{equation}
which differs from the previously~\cite{Biros2005a} derived
\begin{equation}
\left(\mbf{P}_2^{-1} \KKTmat \right)_{\text{BG}}
=
\begin{bmatrix}
 \identity_{\state} & \bzero & \bzero
\\
 \Lyy^T \Rs^{-1} & \Lzz \approxHessian^{-1}  & \bzero
\\
 \Lss \Rs^{-1} & \Lyy \approxHessian^{-1}  &  \identity_{\state}
\end{bmatrix}
.
\label{eq:p2invK_error}
\end{equation}
Both of those results have been verified numerically and algebraically in \ref{app:precond_proof}.

Furthermore, the eigenvalue spectrum of the amended preconditioned system are the same as the previously defined ones
% non preprint
%\begin{equation}
%    \mathcal{E}\left(\mbf{P}^{-1}_4 \KKTmat\right)
%    = \mathcal{E}\left(\mbf{P}^{-1}_2 \KKTmat\right)
%    = \mathcal{E}\left(\left(\mbf{P}_4^{-1} \KKTmat \right)_{\text{BG}}\right)
%    = \mathcal{E}\left(\left(\mbf{P}_2^{-1} \KKTmat \right)_{\text{BG}}\right)
%    = 
%    \left\{ 
%    \left\{ \mathcal{E}(\approxHessian^{-1} \lagrangian_{zz}) \right\}
%    ,  \left\{ \mathbf{\identity_s} \right\}
%    ,  \left\{ \mathbf{\identity_s} \right\}
%    \right\}
%    .
%\end{equation}
%preprint
\begin{equation}
    \mathcal{E}\left(\mbf{P}^{-1}_4 \KKTmat\right)
    = \mathcal{E}\left(\left(\mbf{P}_4^{-1} \KKTmat \right)_{\text{BG}}\right)
    = 
    \left\{ 
    \left\{ \mathcal{E}(\approxHessian^{-1} \lagrangian_{zz}) \right\}
    ,  \left\{ \mathbf{\identity_s} \right\}
    ,  \left\{ \mathbf{\identity_s} \right\}
    \right\}
    .
\end{equation}
\begin{equation}
    \mathcal{E}\left(\mbf{P}^{-1}_2 \KKTmat\right)
    = \mathcal{E}\left(\left(\mbf{P}_2^{-1} \KKTmat \right)_{\text{BG}}\right)
    = 
    \left\{ 
    \left\{ \mathcal{E}(\approxHessian^{-1} \lagrangian_{zz}) \right\}
    ,  \left\{ \mathbf{\identity_s} \right\}
    ,  \left\{ \mathbf{\identity_s} \right\}
    \right\}
    .
\end{equation}
%The eigenvalues have been verified numerically, but further work is required to derive them algebraically.
The eigenvalues have been verified numerically
and algebraically in \ref{app:precond_eig}, which means the results from \cite{Biros2005a} still hold true.

\clearpage

%% file: 9z_precond_eig.tex
\section{Preconditioned system eigenvalues}
\label{app:precond_eig}

To obtain the eigenvalues of the preconditioned system $\mbf{P}^{-1}_4 \KKTmat$, we need to find the roots of $\det \left(\mbf{P}^{-1}_4 \KKTmat - \lambda \identity \right) = 0$
\begin{equation}
\begin{split}
    \det \left(\mbf{P}^{-1}_4 \KKTmat - \lambda \identity \right)
    &=
    \det\left(
    \begin{bmatrix}
\identity_{\state} - \lambda \identity_{\state} 
&
\Rs^{-1}  \Rx
    \left(
    \identity_{\des} - 
    \approxHessian^{-1} \Lzz
    \right)
& \bzero \\
\bzero & \approxHessian^{-1} \Lzz - \lambda \identity_{\des} & \bzero \\
\bzero 
&
  \Rs^{-T} \Lyy
  \left(
  \identity_{\des} - \approxHessian^{-1} \Lzz
  \right)
&
\identity_{\state} - \lambda \identity_{\state} 
\end{bmatrix}
    \right)
\\&=
    \det\left(
    \begin{bmatrix}
    \mbf{M}
& \bzero \\
\mbf{N}
&
\identity_{\state} - \lambda \identity_{\state} 
\end{bmatrix}
    \right),
\end{split}
\end{equation}
where
\begin{equation}
    \mbf{M} =
    \begin{bmatrix}
    \identity_{\state} - \lambda \identity_{\state} 
&
\Rs^{-1}  \Rx
    \left(
    \identity_{\des} - 
    \approxHessian^{-1} \Lzz
    \right) \\
\bzero & \approxHessian^{-1} \Lzz - \lambda \identity_{\des}  \\
    \end{bmatrix}
\quad \text{ and } \quad
\mbf{N} =
\begin{bmatrix}
\bzero 
&
  \Rs^{-T} \Lyy
  \left(
  \identity_{\des} - \approxHessian^{-1} \Lzz
  \right)
\end{bmatrix}
.
\end{equation}
Therefore, the eigenvalues of $\mbf{P}^{-1}_4 \KKTmat$ satisfy the roots of
\begin{equation}
\begin{split}
    \det \left(\mbf{P}^{-1}_4 \KKTmat - \lambda \identity \right)
    &=
    \det\left(\mbf{M}\right)
    \det\left(\identity_{\state} - \lambda \identity_{\state} \right)
    \\&=
    \det\left( \identity_{\state} - \lambda \identity_{\state} \right)
    \det\left( \approxHessian^{-1} \Lzz - \lambda \identity_{\des} \right)
    \det\left( \identity_{\state} - \lambda \identity_{\state} \right) = 0
.
\end{split}
\end{equation}
The preconditioned system $\mbf{P}^{-1}_4 \KKTmat$ therefore has the root 1 with a multiplicity of $2\nstate$, and the other roots corresponding to the eigenvalues of $\approxHessian^{-1} \Lzz$.

We repeat the process for the eigenvalues of $\mbf{P}_2^{-1} \KKTmat$
\begin{equation}
\begin{split}
& \det\left( \mbf{P}_2^{-1} \KKTmat - \lambda \identity \right) 
\\&\quad=
\det\left(
\begin{bmatrix}
% nonpreprint
% \identity_{\state} -  \lambda \identity_{\state}
% -
% \Rs^{-1} \Rx
% \approxHessian^{-1} \Lyy^T
% preprint
\left(\begin{aligned}
 &\identity_{\state} -  \lambda \identity_{\state}
 \\&\quad
 -
 \Rs^{-1} \Rx
 \approxHessian^{-1} \Lyy^T
\end{aligned}\right)
&
% nonpreprint
%  \Rs^{-1} \Rx
% - \Rs^{-1} \Rx \approxHessian^{-1}  
% \left(
% \Lzz 
% -
% \Lyy^T \Rs^{-1} \Rx
% \right)
% preprint
\left(
\begin{aligned}
  &\Rs^{-1} \Rx
 \\&\quad- \Rs^{-1} \Rx \approxHessian^{-1}  
 \left(
 \Lzz 
 -
 \Lyy^T \Rs^{-1} \Rx
 \right)
\end{aligned}
\right)
&
\bzero
\\
 \approxHessian^{-1} \Lyy^T
&
 \approxHessian^{-1} 
 \left(
 \Lzz 
 + \Lyy^T \Rs^{-1} \Rx
 \right)
 - \lambda \identity_{\des}
&
\bzero
\\
 \Rs^{-T}  \Lss
&
 \Rs^{-T} 
 \Lsx
&
 \identity_{\state} -  \lambda \identity_{\state}
\end{bmatrix}
\right)
\\&\quad=
\det\left(
\begin{bmatrix}
\identity_{\state}  &  -\Rs^{-1} \Rx & \bzero \\
\bzero  & \identity_{\des}   & \bzero \\
\bzero  & \bzero & \identity_{\state}
\end{bmatrix}
\begin{bmatrix}
% Row 1
 \identity_{\state} -  \lambda \identity_{\state}
&
 \Rs^{-1} \Rx \left( \identity_{\des} - \lambda \identity_{\des} \right)
&
\bzero
% Row 2
\\
 \approxHessian^{-1} \Lyy^T
&
 \approxHessian^{-1} 
 \left(
 \Lzz 
 + \Lyy^T \Rs^{-1} \Rx
 \right)
 - \lambda \identity_{\des}
&
\bzero
% Row 3
\\
 \Rs^{-T}  \Lss
&
 \Rs^{-T} 
 \Lsx
&
 \identity_{\state} -  \lambda \identity_{\state}
\end{bmatrix}
\right)
\\&\quad=
\det\left(
\begin{bmatrix}
\identity_{\state}  &  -\Rs^{-1} \Rx & \bzero \\
\bzero  & \identity_{\des}   & \bzero \\
\bzero  & \bzero & \identity_{\state}
\end{bmatrix}
\right)
\det
\left(
\begin{bmatrix}
% Row 1
 \mbf{R} &
\bzero
% Row 3
\\
\mbf{S}
&
 \identity_{\state} -  \lambda \identity_{\state}
\end{bmatrix}
\right)
.
\end{split}
\label{eq:p2invK_fact}
\end{equation}
Let,
\begin{equation}
    \mbf{R} =
    \begin{bmatrix}
 \identity_{\state} -  \lambda \identity_{\state}
&
 \Rs^{-1} \Rx \left( \identity_{\des} - \lambda \identity_{\des} \right)
\\
 \approxHessian^{-1} \Lyy^T
&
 \approxHessian^{-1} 
 \left(
 \Lzz 
 + \Lyy^T \Rs^{-1} \Rx
 \right)
 - \lambda \identity_{\des}
 \end{bmatrix}
\quad \text{ and } \quad
    \mbf{S} =
    \begin{bmatrix}
     \Rs^{-T}  \Lss
&
 \Rs^{-T} 
 \Lsx
    \end{bmatrix}
.
\end{equation}
The eigenvalues of $\mbf{P}_2^{-1} \KKTmat$ are the roots of
\begin{equation}
\begin{split}
\det\left( \mbf{P}_2^{-1} \KKTmat - \lambda \identity \right)
&=
\det\left( \identity_{\state} \right)
\det\left( \identity_{\des} \right)
\det\left( \identity_{\state} \right)
\det\left( \mbf{R} \right)
\det\left( \identity_{\state} - \lambda \identity_{\state} \right)
= 0
\end{split}.
\end{equation}
For a matrix
\begin{equation}
    M = \begin{bmatrix}
    A & B \\ C & D
    \end{bmatrix},
\end{equation}
where $A$ is invertible, the determinant of $M$ is given by
\begin{equation}
    \det\left(M\right) = \det\left(A\right) \det\left(D - C A^{-1} B \right)
    \label{eq:schur_complementA}
\end{equation}
Note that if $\lambda = 1$, $\left( \identity_{\state} - \lambda \identity_{\state} \right)$ is not invertible, but $\det\left( \mbf{R} \right) = 0$ nonetheless.
%If $\lambda \neq 1$, we can use Eq. \eqref{eq:schur_complementA} to obtain $\det\left( \mbf{R} \right) = \det\left( \identity_{\state} - \lambda \identity_{\state} \right)
%\det\left( \approxHessian^{-1}  \Lzz  \right)$.
If $\lambda \neq 1$, we can use Eq. \eqref{eq:schur_complementA} to obtain 
\begin{equation}
\begin{split}
    \det\left( \mbf{R} \right) 
    &= 
    \det\left(\identity_{\state} -  \lambda \identity_{\state}\right)
    \det\left(
    \approxHessian^{-1} 
 \left(
 \Lzz 
 + \Lyy^T \Rs^{-1} \Rx
 \right)
 \right.\\&\quad\quad\quad\left.
 - \lambda \identity_{\des}
 -
 \left(\approxHessian^{-1} \Lyy^T \right)
 \left(\identity_{\state} -  \lambda \identity_{\state}\right)^{-1}
 \left( \Rs^{-1} \Rx \left( \identity_{\des} - \lambda \identity_{\des} \right) \right)
\right)
    \\&= 
    \det\left(\identity_{\state} -  \lambda \identity_{\state}\right)
    \det\left(
    \approxHessian^{-1}  \Lzz 
    +
    \approxHessian^{-1}  \Lyy^T \Rs^{-1} \Rx
 \right.\\&\quad\quad\quad\left.
 - \lambda \identity_{\des}
 -
 \left(\approxHessian^{-1} \Lyy^T \right)
 \left( \Rs^{-1} \Rx \right)
\right)
\\&= 
    \det\left(\identity_{\state} -  \lambda \identity_{\state}\right)
    \det\left(
    \approxHessian^{-1}  \Lzz 
 - \lambda \identity_{\des}
\right)
.
\end{split}
\end{equation}
Therefore,
\begin{equation}
\begin{split}
\det\left( \mbf{P}_2^{-1} \KKTmat - \lambda \identity \right)
&=
\det\left( \mbf{R} \right)
\det\left( \identity_{\state} - \lambda \identity_{\state} \right)
\\&=
\det\left( \identity_{\state} - \lambda \identity_{\state} \right)
\det\left( \approxHessian^{-1}  \Lzz  - \lambda \identity_{\des} \right)
\det\left( \identity_{\state} - \lambda \identity_{\state} \right)
= 0
\end{split}
\end{equation}
results in the same roots as the previous preconditioned system.

\clearpage

%% file: 9z_precond_proof.tex
\section{Preconditioned system derivation}
\label{app:precond_proof}
To derive the inverses, we start by checking what the $\mbf{P}_4$ operator can be represented as
\begin{equation} \begin{alignedat}{-1}
\mbf{P}_4
& = \KKTmat_1 \KKTmat_2
&& =
\begin{bmatrix}
\Lss \Rs^{-1} & \bzero & \identity \\
\Lxs \Rs^{-1} & \identity & \Rx^T \Rs^{-T} \\
\identity & \bzero & \bzero
\end{bmatrix}
\begin{bmatrix}
\Rs & \Rx & \bzero \\
\bzero & \approxHessian & \bzero \\
\bzero &
\left(  
    \Lsx-\Lss \Rs^{-1} \Rx 
\right)
& \Rs^T
\end{bmatrix}
\\
&&& =
\begin{bmatrix}
\Lss
&
\Lsx
&
\Rs^{T}
\\
\Lxs
&
\Lxs \Rs^{-1} \Rx + \approxHessian + \Rx^T \Rs^{-T}
\left(  
    \Lsx-\Lss \Rs^{-1} \Rx 
\right)
&
\Rx^{T}
\\
\Rs & \Rx & \bzero
\end{bmatrix}
\\
&&& =
\begin{bmatrix}
\Lss
&
\Lsx
&
\Rs^{T}
\\
\Lxs
&
\approxHessian - \Lzz + \Lxx
&
\Rx^{T}
\\
\Rs & \Rx & \bzero
\end{bmatrix}
.
\end{alignedat} \end{equation}
Therefore, if $\approxHessian = \Lzz $, we recover the KKT matrix as expected.

Let's find $\KKTmat_1^{-1}$
\begin{equation}
\begin{bmatrix}
\Lss \Rs^{-1} & \bzero & \identity \\
\Lxs \Rs^{-1} & \identity & \Rx^T \Rs^{-T} \\
\identity & \bzero & \bzero
\end{bmatrix}
\begin{bmatrix}
\KoneInv_{11} & \KoneInv_{12} & \KoneInv_{13} \\
\KoneInv_{21} & \KoneInv_{22} & \KoneInv_{23} \\
\KoneInv_{31} & \KoneInv_{32} & \KoneInv_{33}
\end{bmatrix}
=
\begin{bmatrix}
\identity_{\state} & \bzero & \bzero \\
\bzero & \identity_{\des} & \bzero \\
\bzero & \bzero & \identity_{\state}
\end{bmatrix}
.
\end{equation}

%%% \begin{equation}
%%% \begin{bmatrix}
%%% \Lss \Rs^{-1} & \bzero & \identity \\
%%% \Lxs \Rs^{-1} & \identity & \Rx^T \Rs^{-T} \\
%%% \identity & \bzero & \bzero
%%% \end{bmatrix}
%%% \begin{bmatrix}
%%% \KoneInv_{11} \\
%%% \KoneInv_{21} \\
%%% \KoneInv_{31} 
%%% \end{bmatrix}
%%% =
%%% \begin{bmatrix}
%%% \identity_{\state} \\
%%% \bzero        \\
%%% \bzero        
%%% \end{bmatrix}
%%% ,\quad
%%% \begin{bmatrix}
%%% \Lss \Rs^{-1} & \bzero & \identity \\
%%% \Lxs \Rs^{-1} & \identity & \Rx^T \Rs^{-T} \\
%%% \identity & \bzero & \bzero
%%% \end{bmatrix}
%%% \begin{bmatrix}
%%% \KoneInv_{12} \\
%%% \KoneInv_{22} \\
%%% \KoneInv_{32} 
%%% \end{bmatrix}
%%% =
%%% \begin{bmatrix}
%%% \bzero        \\
%%% \identity_{\des} \\
%%% \bzero        
%%% \end{bmatrix}
%%% ,\quad
%%% \begin{bmatrix}
%%% \Lss \Rs^{-1} & \bzero & \identity \\
%%% \Lxs \Rs^{-1} & \identity & \Rx^T \Rs^{-T} \\
%%% \identity & \bzero & \bzero
%%% \end{bmatrix}
%%% \begin{bmatrix}
%%% \KoneInv_{13} \\
%%% \KoneInv_{23} \\
%%% \KoneInv_{33} 
%%% \end{bmatrix}
%%% =
%%% \begin{bmatrix}
%%% \bzero        \\
%%% \bzero        \\
%%% \identity_{\state}
%%% \end{bmatrix}
%%% \end{equation}

%% First column
The first column forms
\begin{equation} \begin{alignedat}{-1}
% First row
{ \Lss \Rs^{-1} \KoneInv_{11} } &
{} &
+ { \identity_{\state} \KoneInv_{31} } &
&=&&\identity_{\state}
,\\
% Second row
{ \Lxs \Rs^{-1} \KoneInv_{11} } &
+ { \identity_{\des} \KoneInv_{21} } &
+ { \Rx^T \Rs^{-T} \KoneInv_{31} } &
&=&& \bzero
,\\
% Third row
{ \identity_{\state} \KoneInv_{11} } &
{ } &
{ } &
&=&& \bzero
,
\end{alignedat} \end{equation}
which gives 
\begin{equation} \begin{alignedat}{2}
\KoneInv_{11} 
&&& = \bzero
,\\
\KoneInv_{31}
&  = \left(\identity_{\state} - \Lss \Rs^{-1} \mbf{T}_{11} \right)
&& = \identity_{\state}
,\\
\KoneInv_{21} 
&  = 
 - \Lxs \Rs^{-1} \KoneInv_{11} 
 - \Rx^T \Rs^{-T} \KoneInv_{31} 
&& = - \Rx^T \Rs^{-T}
.
\end{alignedat} \end{equation}
%
%% Second column
The second column forms
\begin{equation} \begin{alignedat}{-1}
% First row
{ \Lss \Rs^{-1} \KoneInv_{12} } &
{} &
+ { \identity_{\state} \KoneInv_{32} } &
&=&& \bzero
,\\
% Second row
{ \Lxs \Rs^{-1} \KoneInv_{12} } &
+ { \identity_{\des} \KoneInv_{22} } &
+ { \Rx^T \Rs^{-T} \KoneInv_{32} } &
&=&& \identity_{\des}
,\\
% Third row
{ \identity_{\state} \KoneInv_{12} } &
{ } &
{ } &
&=&& \bzero
,
\end{alignedat} \end{equation}
which gives 
\begin{equation} \begin{alignedat}{2}
\KoneInv_{12} 
&&& = \bzero
,\\
\KoneInv_{32}
&  = - \Lss \Rs^{-1} \KoneInv_{12}
&& = \bzero
,\\
\KoneInv_{22} 
&  = \identity_{\des} 
 - \Lxs \Rs^{-1} \KoneInv_{12} 
 - \Rx^T \Rs^{-T} \KoneInv_{32} 
&& = \identity_{\des}
.
\end{alignedat} \end{equation}
%% Third column
%
The third column forms
\begin{equation} \begin{alignedat}{-1}
% First row
{ \Lss \Rs^{-1} \KoneInv_{13} } &
{} &
+ { \identity_{\state} \KoneInv_{33} } &
&=&& \bzero
,\\
% Second row
{ \Lxs \Rs^{-1} \KoneInv_{13} } &
+ { \identity_{\des} \KoneInv_{23} } &
+ { \Rx^T \Rs^{-T} \KoneInv_{33} } &
&=&& \bzero
,\\
% Third row
{ \identity_{\state} \KoneInv_{13} } &
{ } &
{ } &
&=&& \identity_{\state}
,
\end{alignedat} \end{equation}
which gives 
\begin{equation} \begin{alignedat}{2}
\KoneInv_{13} 
&&& = \identity_{\state}
,\\
\KoneInv_{33}
&  = - \Lss \Rs^{-1} \KoneInv_{13} 
&& = - \Lss \Rs^{-1} 
,\\
\KoneInv_{23} 
&  = 
 - \Lxs \Rs^{-1} \KoneInv_{13} 
 - \Rx^T \Rs^{-T} \KoneInv_{33} 
\\& = 
%&& = % non preprint alignment
- \Lxs \Rs^{-1} 
 + \Rx^T \Rs^{-T} \Lss \Rs^{-1} 
 = -\Lyy^T \Rs^{-1} 
.
\end{alignedat}\end{equation}

Repeat the whole process for $\KKTmat_2$.
%% First column
The first column gives 
\begin{equation} \begin{alignedat}{2}
\KtwoInv_{21} 
& 
&& = \bzero
,
\\
\KtwoInv_{11} 
&  = \Rs^{-1} \left(\identity_{\state} - \Rx \KtwoInv_{21} \right)
&& = \Rs^{-1}
,
\\
\KtwoInv_{31}
&  = -\Rs^{-T}
\Lyy
\KtwoInv_{21}
&& = \bzero
.
\end{alignedat}\end{equation}
The second column results in
\begin{equation} \begin{alignedat}{2}
\KtwoInv_{22} 
&  = \approxHessian^{-1} \identity_{\des}
&& = \approxHessian^{-1}
,\\
\KtwoInv_{12} 
&  = -\Rs^{-1} \Rx \KtwoInv_{22} 
&& = -\Rs^{-1} \Rx \approxHessian^{-1}
,\\
\KtwoInv_{32}
&  = -\Rs^{-T}
\Lyy
%\left(\Lsx - \Lss \Rs^{-1} \Rx \right)
\KtwoInv_{22}
&& = -\Rs^{-T}
\Lyy
%\left(\Lsx - \Lss \Rs^{-1} \Rx \right) 
\approxHessian^{-1}
.
\end{alignedat}\end{equation}
The third column results in
\begin{equation} \begin{alignedat}{2}
\KtwoInv_{23} 
&  
&& = \bzero
,\\
\KtwoInv_{13} 
&  = -\Rs^{-1} \Rx \KtwoInv_{23} 
&& = \bzero
,\\
\KtwoInv_{33}
&  = \Rs^{-T} 
\left( 
    \identity_{\state}
    - 
    \Lyy
    %\left(\Lsx - \Lss \Rs^{-1} \Rx \right)
    \KtwoInv_{23}
\right)
&& = \Rs^{-T} 
.
\end{alignedat}\end{equation}

Finally, the inverse of $\KKTmat_1$ is given by
\begin{equation}
\KKTmat_1^{-1} = 
\begin{bmatrix}
 \bzero & \bzero & \identity_{\state}
\\
 -\Rx^T \Rs^{-T}
 &
 \identity_{\des}
 &
 -\Lyy^T \Rs^{-1} 
 \\
 \identity_{\state}
 &
 \bzero
 &
 -\Lss \Rs^{-1}
\end{bmatrix}
,
\end{equation}
and the inverse of $\KKTmat_2$ is given by
\begin{equation}
\KKTmat_2^{-1} = 
\begin{bmatrix}
 \Rs^{-1} 
 & -\Rs^{-1} \Rx \approxHessian^{-1}
 & \bzero
\\
 \bzero & \approxHessian^{-1} & \bzero
 \\
 \bzero
 &
  -\Rs^{-T} 
  \Lyy
  %\left(\Lsx - \Lss \Rs^{-1} \Rx \right) 
  \approxHessian^{-1}
 &
 \Rs^{-T}
\end{bmatrix}
.
\end{equation}
Finally, the inverse of the $\mbf{P}_4$ preconditioner is given by
\begin{equation}
\mbf{P}_4^{-1}
= \left( \KKTmat_1 \KKTmat_2 \right)^{-1}
= \KKTmat_2^{-1} \KKTmat_1^{-1}.
\end{equation}
\begin{equation}
\begin{split}
\PfourInv_{11} &= 
    \Rs^{-1}\Rx  \approxHessian^{-1}
    \Rx^T \Rs^{-T}
,\\ 
\PfourInv_{12} &= 
    -\Rs^{-1}\Rx  \approxHessian^{-1}
,\\ 
\PfourInv_{13} &= 
    \Rs^{-1} 
    + \Rs^{-1}\Rx  \approxHessian^{-1}
        \Lyy^T \Rs^{-1} 
,\\
\PfourInv_{21} &= 
    -\approxHessian^{-1} \Rx^T \Rs^{-T}
,\\ 
\PfourInv_{22} &=  \approxHessian^{-1}
,\\
\PfourInv_{23} &=  -\approxHessian^{-1}
        \Lyy^T \Rs^{-1} 
,\\
\PfourInv_{31} &= 
  \Rs^{-T}
  \Lyy
  %\left(\Lsx - \Lss \Rs^{-1} \Rx \right)
  \approxHessian^{-1}
  \Rx^T \Rs^{-T}
  -\Rs^{-T}
,\\ \PfourInv_{32} &= 
  - \Rs^{-T} 
  \Lyy
  %\left(\Lsx - \Lss \Rs^{-1} \Rx \right)
  \approxHessian^{-1} 
,\\ \PfourInv_{33} &= 
  \Rs^{-T} \Lyy
  \approxHessian^{-1}
  \Lyy^T \Rs^{-1} 
  -\Rs^{-T}
  \Lss\Rs^{-1}
.
\end{split}
\end{equation}

%Interestingly, the eigenvalues still only depend on the diagonal for the given block pattern.
%\begin{equation}
%\begin{vmatrix}
%\mbf{C}_{11} - \lambda \identity_{\state} & \mbf{C}_{12} & \bzero \\
%\bzero & \mbf{C}_{22} - \lambda \identity_{\des} & \bzero \\
%\bzero & \mbf{C}_{32} & \mbf{C}_{33} - \lambda \identity_{\state}
%\end{vmatrix}
%=
%\left| \mbf{C}_{11} - \lambda \identity_{\state} \right|
%\left( \mbf{C}_{22} - \lambda \identity_{\des} \right| 
%\left( \mbf{C}_{33} - \lambda \identity_{\state} \right|
%\end{equation}
Multiplying the inverse of the preconditioner $\mbf{P}_4^{-1}$ with the KKT matrix $\KKTmat$ gives the following preconditioned system
\begin{equation}
    \mbf{P}_4^{-1} \KKTmat = \mbf{C}
\end{equation}

\begin{equation} \begin{alignedat}{-1}
\mbf{C}_{11}
&= 
    \Rs^{-1}\Rx  \approxHessian^{-1}
    \Rx^T \Rs^{-T}
    \Lss
    - \Rs^{-1}\Rx  \approxHessian^{-1}
    \Lxs
\\&\quad\quad +
    \left(\Rs^{-1} 
    + \Rs^{-1}\Rx  \approxHessian^{-1}
    \Lyy^{T} \Rs^{-1} 
    \right)
    \Rs
\\&=
    \Rs^{-1}\Rx  \approxHessian^{-1}
    \Rx^T \Rs^{-T}
    \Lss
-
    \Rs^{-1}\Rx  \approxHessian^{-1}
    \Lxs
+
    \identity_{\state}
    + \Rs^{-1}\Rx  \approxHessian^{-1}
    \Lyy^{T}
\\&=
    \Rs^{-1}\Rx  \approxHessian^{-1}
    \left(
        \Rx^T \Rs^{-T} \Lss
        - \Lxs
    \right)
    + \identity_{\state}
    + \Rs^{-1}\Rx  \approxHessian^{-1}
    \Lyy^{T}
\\&=
    -\Rs^{-1}\Rx  \approxHessian^{-1}
    \Lyy^T
    + \identity_{\state}
    + \Rs^{-1}\Rx  \approxHessian^{-1}
    \Lyy^{T}
\\&= \identity_{\state}
\\
\mbf{C}_{12}
&= 
    \Rs^{-1}\Rx  \approxHessian^{-1}
    \Rx^T \Rs^{-T}
    \Lsx
-
    \Rs^{-1}\Rx  \approxHessian^{-1}
    \Lxx
\\&\quad\quad +
    \left(\Rs^{-1} 
    + \Rs^{-1}\Rx  \approxHessian^{-1}
    \Lyy^{T} \Rs^{-1} 
    \right)
    \Rx
\\&= 
    \Rs^{-1}\Rx  \approxHessian^{-1}
    \Rx^T \Rs^{-T}
    \Lsx
-
    \Rs^{-1}\Rx  \approxHessian^{-1}
    \Lxx
\\&\quad\quad +
    \Rs^{-1} 
    \Rx
    + \Rs^{-1}\Rx  \approxHessian^{-1}
    \Lyy^{T} \Rs^{-1} 
    \Rx
\\&= 
    \Rs^{-1}\Rx  \approxHessian^{-1}
    \left(
        \Rx^T \Rs^{-T} \Lsx
        - \Lxx
        + \Lyy^{T} \Rs^{-1}  \Rx
    \right)
    + \Rs^{-1}  \Rx
\\&= 
    \Rs^{-1}\Rx  \approxHessian^{-1}
    \left(-\Lzz \right)
    + \Rs^{-1}  \Rx
\\&= 
    \Rs^{-1}  \Rx
    \left(
    \identity_{\des} - 
    \approxHessian^{-1} \Lzz
    \right)
\\
\mbf{C}_{13}
&= 
    \Rs^{-1}\Rx  \approxHessian^{-1}
    \Rx^T \Rs^{-T}
    \Rs^T
-
    \Rs^{-1}\Rx  \approxHessian^{-1}
    \Rx^T
\\&= 
    \Rs^{-1}\Rx  \approxHessian^{-1}
    \Rx^T 
-
    \Rs^{-1}\Rx  \approxHessian^{-1}
    \Rx^T
\\&=  \bzero
\end{alignedat}\end{equation}
\begin{equation} \begin{alignedat}{-1}
\mbf{C}_{21}
&= 
    -\approxHessian^{-1} \Rx^T \Rs^{-T}
    \Lss
+ 
\approxHessian^{-1}
    \Lxs
- \approxHessian^{-1}
        \Lyy^T \Rs^{-1} 
    \Rs
\\&= 
\approxHessian^{-1}
        \Lyy^T 
- \approxHessian^{-1}
        \Lyy^T 
\\&= \bzero
\\
\mbf{C}_{22}
&= 
    -\approxHessian^{-1} \Rx^T \Rs^{-T}
    \Lsx
+ 
\approxHessian^{-1}
    \Lxx
- \approxHessian^{-1}
        \Lyy^T \Rs^{-1} 
    \Rx
\\&= 
\approxHessian^{-1}
\left(
    - \Rx^T \Rs^{-T} \Lsx
    +  \Lxx
    -  \Lyy^T \Rs^{-1}  \Rx
\right)
\\&= 
\approxHessian^{-1} \lagrangian_{z}
\\
\mbf{C}_{23}
&= 
    -\approxHessian^{-1} \Rx^T \Rs^{-T}
    \Rs^T
+ 
\approxHessian^{-1}
    \Rx^T
\\&= \bzero
\end{alignedat}\end{equation}
\begin{equation} \begin{alignedat}{-1}
\mbf{C}_{31}
&= 
  \left(
      \Rs^{-T} \Lyy \approxHessian^{-1}
          \Rx^T \Rs^{-T}
      + \Rs^{-T}
  \right)
  \Lss
- \Rs^{-T}  \Lyy \approxHessian^{-1} 
  \Lxs
\\&\quad\quad+
\left(
  \Rs^{-T} \Lyy
  \approxHessian^{-1}
  \Lyy^T \Rs^{-1} 
  -\Rs^{-T}
  \Lss\Rs^{-1}
 \right)
 \Rs
\\&= 
  \left(
      \Rs^{-T} \Lyy \approxHessian^{-1}
          \Rx^T 
      + \identity_{\state}
  \right)
  \Rs^{-T}
  \Lss
- \Rs^{-T}  \Lyy \approxHessian^{-1} 
  \Lxs
\\&\quad\quad+
  \Rs^{-T} \Lyy
  \approxHessian^{-1}
  \Lyy^T
  -\Rs^{-T}
  \Lss
\\&= 
  \Rs^{-T} \Lyy \approxHessian^{-1}
  \Rx^T  \Rs^{-T} \Lss
- \Rs^{-T}  \Lyy \approxHessian^{-1} 
  \Lxs
+
  \Rs^{-T} \Lyy
  \approxHessian^{-1}
  \Lyy^T
\\&= 
  \Rs^{-T} \Lyy \approxHessian^{-1}
  \left(
      \Rx^T  \Rs^{-T}
      \Lss -  \Lxs
   \right)
+
  \Rs^{-T} \Lyy
  \approxHessian^{-1}
  \Lyy^T
\\&= 
  -\Rs^{-T} \Lyy \approxHessian^{-1}
  \Lyy^T
+
  \Rs^{-T} \Lyy
  \approxHessian^{-1}
  \Lyy^T
\\&= \bzero
\\
\mbf{C}_{32}
&= 
  \left(
      \Rs^{-T} \Lyy \approxHessian^{-1}
          \Rx^T \Rs^{-T}
      + \Rs^{-T}
  \right)
  \Lsx
- \Rs^{-T}  \Lyy \approxHessian^{-1} 
  \Lxx
\\&\quad\quad+
\left(
  \Rs^{-T} \Lyy
  \approxHessian^{-1}
  \Lyy^T \Rs^{-1} 
  -\Rs^{-T}
  \Lss\Rs^{-1}
 \right)
 \Rx
\\&= 
  \Rs^{-T} \Lyy \approxHessian^{-1}
  \Rx^T \Rs^{-T}
  \Lsx
  + \Rs^{-T} \Lsx
- \Rs^{-T}  \Lyy \approxHessian^{-1} 
  \Lxx
\\&\quad\quad+
  \Rs^{-T} \Lyy
  \approxHessian^{-1}
  \Lyy^T \Rs^{-1} 
 \Rx
  -\Rs^{-T}
  \Lss\Rs^{-1}
 \Rx
\\&= 
  \Rs^{-T} \Lyy \approxHessian^{-1}
  \left(
    \Rx^T \Rs^{-T} \Lsx
    - \Lxx
    + \Lyy^T \Rs^{-1}  \Rx
  \right)
\\&\quad\quad
  + \Rs^{-T} \Lsx
  -\Rs^{-T} \Lss\Rs^{-1}
 \Rx
\\&= 
  \Rs^{-T} \Lyy \approxHessian^{-1}
  \left(
    \Rx^T \Rs^{-T} \Lsx
    - \Lxx
    + \Lyy^T \Rs^{-1}  \Rx
  \right)
  + \Rs^{-T} \Lyy
\\&= 
  \Rs^{-T} \Lyy \approxHessian^{-1}
  \left(
      -\Lzz
  \right)
  + \Rs^{-T} \Lyy
\\&= 
  -\Rs^{-T} \Lyy
  \approxHessian^{-1} \Lzz
  + \Rs^{-T} \Lyy
\\&= 
  \Rs^{-T} \Lyy
  \left(
  \identity_{\des} - \approxHessian^{-1} \Lzz
  \right)
\\
\mbf{C}_{33}
&= 
  \left(
      \Rs^{-T} \Lyy \approxHessian^{-1}
          \Rx^T \Rs^{-T}
      + \Rs^{-T}
  \right)
  \Rs^T
- \Rs^{-T}  \Lyy \approxHessian^{-1} 
  \Rx^T
\\&= 
      \Rs^{-T} \Lyy \approxHessian^{-1}
          \Rx^T 
      + \identity_{\state}
- \Rs^{-T}  \Lyy \approxHessian^{-1} 
  \Rx^T
\\&=  \identity_{\state}
\end{alignedat}\end{equation}

In order to derive $\mbf{P}_2^{-1}$, it is easier to form $\mbf{P}_2$
\begin{equation}
\mbf{P}_2 =
\begin{bmatrix}
%\Lss & \Lsx
\bzero & \bzero
&
\Rs^{T}
\\
%\Lxs
\bzero
&
\approxHessian
&
\Rx^{T}
\\
\Rs & \Rx & \bzero
\end{bmatrix}
,
\end{equation}
and repeat the above process to form its inverse
\begin{equation}
\mbf{P}_2^{-1} = 
\begin{bmatrix}
 \Rs^{-1} \Rx \Lzz^{-1} \Rx^{T} \Rs^{-T}
 &
 - \Rs^{-1} \Rx \Lzz^{-1} 
 &
 \Rs^{-1}
 \\
 - \Lzz^{-1} \Rx^{T} \Rs^{-T} 
 &
 \Lzz^{-1} 
 &
 \bzero
 \\
 \Rs^{-T}
 &
 \bzero
 &
 \bzero
\end{bmatrix}
.
\end{equation}
Finally, apply it onto the KKT matrix to obtain \eqref{eq:p2invK} repeated below
\begin{equation}
\mbf{P}_2^{-1} \KKTmat =
\begin{bmatrix}
% V11
 \identity_{\state}
 -
 \Rs^{-1} \Rx
 \approxHessian^{-1} \Lyy^T
&
% V12
 \Rs^{-1} \Rx
 - \Rs^{-1} \Rx \approxHessian^{-1}  
 \left(
 \Lzz 
 -
 \Lyy^T \Rs^{-1} \Rx
 \right)
&
% V13
\bzero
\\
% V21
 \approxHessian^{-1} \Lyy^T
&
% V22
 \approxHessian^{-1} 
 \left(
 \Lzz 
 +\Lyy^T \Rs^{-1} \Rx
 \right)
&
% V23
\bzero
\\
% V31
 \Rs^{-T} \Lss
&
% V32
 \Rs^{-T} 
\Lsx
&
% 33
\identity_{\state}
\end{bmatrix}
.
\end{equation}

%   \begin{equation}
%   \mbf{P}_2^{-1} \KKTmat =
%   \begin{bmatrix}
%   % V11
%    \identity_{\state}
%    -
%    \Rs^{-1} \Rx
%    \approxHessian^{-1} \Lyy^T
%   &
%   % V12
%    \Rs^{-1} \Rx
%    \left(
%        \identity_{\des}
%        - \approxHessian^{-1}  \Lzz
%    \right)
%    -
%    \Rs^{-1} \Rx
%    \approxHessian^{-1}  \Lyy^T
%    \Rs^{-1} \Rx
%   &
%   % V13
%   \bzero
%   \\
%   % V21
%    \approxHessian^{-1} \Lyy^T
%   &
%   % V22
%    \approxHessian^{-1} 
%    \left(
%    \Lzz 
%    +\Lyy^T \Rs^{-1} \Rx
%    \right)
%   &
%   % V23
%   \bzero
%   \\
%   % V31
%    \Rs^{-T} \Lss
%   &
%   % V32
%    \Rs^{-T} 
%   \Lsx
%   &
%   % 33
%   \identity_{\state}
%   \end{bmatrix}
%   .
%   \end{equation}

%The eigenvalue spectrum of the resulting preconditioned system is therefore not obvious unless further simplifications are made.
%However, numerical experiments have shown that the spectrum of matrix \eqref{eq:p2invK} does match \eqref{eq:p4invK}.

\clearpage

%% file: 9z_flops.tex
\section{Residual flops}
\label{app:residual_flops}
\subsection{Domain contribution}
\label{app:domain_residual_flops}
The number of flops is first counted for a single quadrature point, and then multiplied by the total number of quadrature points within an element.
\begin{itemize}
    \item Interpolating the solution requires $(2s\nb - s)$ flops.
    \item Interpolating the metric gradient to form the metric Jacobian requires $(2d^2\nb - d^2)$ flops.
    \item Forming the cofactor matrix of the Jacobian is approximately $(d^d)$ for $d=1,2,3$ flops.
    \item Evaluating the analytical convective flux requires around $(8d+7)$ flops,
    \item and converting them into the refence domain requires $(2d^2s - ds)$ flops.
    \item Adding its contribution to the $(\nb s)$ residuals by dotting the appropriate basis gradient with the flux requires $(2d s\nb)$.
    \item Finally, multiplying by the quadrature weight requires 1 flop.
\end{itemize}
The above totals
%$(2s\nb - s) + (2d^2\nb - d^2) + d^d + 8d + 7 + 2dds - ds + 2 sd\nb + 1$
%$(2s + 2d^2 + 2ds) \nb + (2dd - d - 1) s + (8 - d) d + d^d + 8
$(4dd + 6d + 4) \nb + d^d + 2ddd + 2dd + 5d + 6$
flops.
In 2D and 3D, we get $32 \nb + 44$ and $58 \nb + 120$ flops per quadrature point.
Multiplying by $\nb$ for its repetitions and adding $\nb$ for the sum. 
We total 
$32\nb^2 + 45\nb$ 
and 
$58\nb^2 + 121\nb$, 
when we add up the $\nb$ quadrature loop for 2D and 3D respectively.

\subsection{Surface contribution}
\label{app:surface_residual_flops}

The number of flops is first counted for a single quadrature point, and then multiplied by the total number of face quadrature points.
\begin{itemize}
    \item Interpolating the solution requires $(2s\nb - s)$ flops.
    \item Interpolating the metric gradient to form the metric Jacobian requires $(2d^2\nb - d^2)$ flops.
    \item Forming the cofactor matrix of the Jacobian  approximately requires $(d^d)$ for $d=1,2,3$ flops.
    \item Contracting the cofactor matrix with unit normal costs $(2d^2-d)$ flops.
    \item Evaluating contravariant numerical Roe flux with entropy fix requires $(46d+130)$ flops.
    \item Adding its contribution to the $2s \nb$ residuals by multiplying the basis values with flux costs $2(2s \nb)$ flops.
    \item Multiplying by quadrature weight for both residual integrals costs 2 flops.
\end{itemize}
The above totals 
%$(2s\nb - s) + (2d^2\nb - d^2) + (d^d) + (2d^2-d) + (46d+130) + 2(2s \nb) + 2$
%$(2s + 2dd + 4s) \nb - s - d^2 + (d^d) + (2d^2-d) + (46d+130) + 2$
%$(2d+4 + 2dd + 4d+8) \nb - d - 2 - d^2 + d^d + 2d^2 - d + 46d + 130 + 2$
$(2dd + 6d + 12) \nb + d^d + dd + 44d + 130$.
In 2D and 3D, we get 
$32 \nb + 226$
and 
$48 \nb + 298$
flops per quadrature point.
Multiplying by $n_f$ quadrature points and adding $n_f$ for the sum, we total 
$32\nb n_f + 227 n_f $
and 
$48 \nb n_f + 299 n_f$.

\subsection{Total residual flops}
A 1/2/3D structured tensor-product mesh will have 1/2/3 times as many faces as there are cells. Therefore, if we consider the cost of assembling the residual for 1 cell, we need to add up the cost of the volume integral plus 1/2/3 times the cost of the face. Another way to see it, is that we compute the face contributions it for all 2/4/6 faces, but divide the work by 2 since only one of the neighbouring cell is responsible for the face computation.

The total flops for a cell's computation in 2D is therefore 
$32\nb^2 + 45\nb + 2(32\nb n_f + 227 n_f)$ flops 
or 
$32(p+1)^4 + 64 (p+1)^{3} + 45(p+1)^2 + 454 (p+1)$
.
In 3D, the total flops is
$58\nb^2 + 121\nb+ 3(48 \nb n_f + 299 n_f)$
or 
$58(p+1)^{6} + 144 (p+1)^{5} + 121 (p+1)^{3} + 897 (p+1)^{2}$

%% file: 00000_main.bbl
\begin{thebibliography}{10}
\expandafter\ifx\csname url\endcsname\relax
  \def\url#1{\texttt{#1}}\fi
\expandafter\ifx\csname urlprefix\endcsname\relax\def\urlprefix{URL }\fi
\expandafter\ifx\csname href\endcsname\relax
  \def\href#1#2{#2} \def\path#1{#1}\fi

\bibitem{Reuther1996}
J.~Reuther, A.~Jameson, J.~Farmer, L.~Martinelli, D.~Saunders, Aerodynamic
  shape optimization of complex aircraft configurations via an adjoint
  formulation, in: 34th Aerospace Sciences Meeting and Exhibit, American
  Institute of Aeronautics and Astronautics, 1996.
\newblock \href {http://dx.doi.org/10.2514/6.1996-94}
  {\path{doi:10.2514/6.1996-94}}.

\bibitem{Palacios2015}
F.~Palacios, T.~D. Economon, J.~J. Alonso, Large-scale aircraft design using
  {SU}2, in: 53rd {AIAA} Aerospace Sciences Meeting, American Institute of
  Aeronautics and Astronautics, 2015.
\newblock \href {http://dx.doi.org/10.2514/6.2015-1946}
  {\path{doi:10.2514/6.2015-1946}}.

\bibitem{Gagnon2016}
H.~Gagnon, D.~W. Zingg, Euler-equation-based drag minimization of
  unconventional aircraft configurations, Journal of Aircraft 53~(5) (2016)
  1361--1371.
\newblock \href {http://dx.doi.org/10.2514/1.c033591}
  {\path{doi:10.2514/1.c033591}}.

\bibitem{Chen2016}
S.~Chen, Z.~Lyu, G.~K.~W. Kenway, J.~R. R.~A. Martins, Aerodynamic shape
  optimization of common research model wing{\textendash}body{\textendash}tail
  configuration, Journal of Aircraft 53~(1) (2016) 276--293.
\newblock \href {http://dx.doi.org/10.2514/1.c033328}
  {\path{doi:10.2514/1.c033328}}.

\bibitem{Luo2010}
J.~Luo, J.~Xiong, F.~Liu, I.~McBean, Three-dimensional aerodynamic design
  optimization of a turbine blade by using an adjoint method, Journal of
  Turbomachinery 133~(1).
\newblock \href {http://dx.doi.org/10.1115/1.4001166}
  {\path{doi:10.1115/1.4001166}}.

\bibitem{Walther2015}
B.~Walther, S.~Nadarajah, Optimum shape design for multirow turbomachinery
  configurations using a discrete adjoint approach and an efficient radial
  basis function deformation scheme for complex multiblock grids, ASME Journal
  of Turbomachinery 137~(8) (2015) 081006.
\newblock \href {http://dx.doi.org/10.1115/1.4029550}
  {\path{doi:10.1115/1.4029550}}.

\bibitem{Jameson1988}
A.~Jameson, Aerodynamic design via control theory, Journal of Scientific
  Computing 3~(3) (1988) 233--260.
\newblock \href {http://dx.doi.org/10.1007/BF01061285}
  {\path{doi:10.1007/BF01061285}}.

\bibitem{Brown2017}
D.~A. Brown, S.~Nadarajah, Inexactly constrained discrete adjoint approach for
  steepest descent-based optimization algorithms, Numerical Algorithms 76~(2)
  (2017) 1--18.
\newblock \href {http://dx.doi.org/10.1007/s11075-017-0409-7}
  {\path{doi:10.1007/s11075-017-0409-7}}.

\bibitem{Sherman1996}
L.~L. Sherman, I.~I.~I. Arthur C.~Taylor, L.~L. Green, P.~A. Newman, G.~W. Hou,
  V.~M. Korivi,
  \href{http://www.sciencedirect.com/science/article/pii/S0021999196902521}{First-
  and second-order aerodynamic sensitivity derivatives via automatic
  differentiation with incremental iterative methods}, Journal of Computational
  Physics 129~(2) (1996) 307--331.
\newblock \href {http://dx.doi.org/10.1006/jcph.1996.0252}
  {\path{doi:10.1006/jcph.1996.0252}}.
\newline\urlprefix\url{http://www.sciencedirect.com/science/article/pii/S0021999196902521}

\bibitem{Ghate2007}
D.~Ghate, M.~Giles, Efficient {Hessian} calculation using automatic
  differentiation, in: 25th AIAA Applied Aerodynamics Conference, 2007, paper
  No. AIAA 2007-4059.
\newblock \href {http://dx.doi.org/10.2514/6.2007-4059}
  {\path{doi:10.2514/6.2007-4059}}.

\bibitem{Papadimitriou2008}
D.~I. Papadimitriou, K.~C. Giannakoglou, Direct, adjoint and mixed approaches
  for the computation of {Hessian} in airfoil design problems, International
  Journal for Numerical Methods in Fluids 56~(10) (2008) 1929--1943.
\newblock \href {http://dx.doi.org/10.1002/fld.1584}
  {\path{doi:10.1002/fld.1584}}.

\bibitem{ShiDong2018}
D.~Shi-Dong, S.~Nadarajah, Approximate {Hessian} for accelerated convergence of
  aerodynamic shape optimization problems in an adjoint-based framework,
  Computers \& Fluids\href {http://dx.doi.org/10.1016/j.compfluid.2018.04.019}
  {\path{doi:10.1016/j.compfluid.2018.04.019}}.

\bibitem{Arian1999}
E.~Arian, S.~Ta'asan,
  \href{http://www.sciencedirect.com/science/article/pii/S0045793098000607}{Analysis
  of the {Hessian} for aerodynamic optimization: inviscid flow}, Computers \&
  Fluids 28~(7) (1999) 853--877.
\newblock \href {http://dx.doi.org/10.1016/s0045-7930(98)00060-7}
  {\path{doi:10.1016/s0045-7930(98)00060-7}}.
\newline\urlprefix\url{http://www.sciencedirect.com/science/article/pii/S0045793098000607}

\bibitem{Schmidt2010}
S.~Schmidt,
  \href{http://ubt.opus.hbz-nrw.de/frontdoor.php?source_opus=569&la=en}{Efficient
  large scale aerodynamic design based on shape calculus}, Ph.D. thesis,
  University of Trier, Germany (Apr. 2010).
\newline\urlprefix\url{http://ubt.opus.hbz-nrw.de/frontdoor.php?source_opus=569&la=en}

\bibitem{Heinkenschloss1999}
M.~Heinkenschloss, L.~N. Vicente, An interface optimization and application for
  the numerical solution of optimal control problems, {ACM} Transactions on
  Mathematical Software 25~(2) (1999) 157--190.
\newblock \href {http://dx.doi.org/10.1145/317275.317278}
  {\path{doi:10.1145/317275.317278}}.

\bibitem{Akcelik2002}
V.~Akcelik, G.~Biros, O.~Ghattas, Parallel multiscale gauss-newton-krylov
  methods for inverse wave propagation, in: {ACM}/{IEEE} {SC} 2002 Conference
  ({SC}{\textsf{'}}02), {IEEE}, 2002.
\newblock \href {http://dx.doi.org/10.1109/sc.2002.10002}
  {\path{doi:10.1109/sc.2002.10002}}.

\bibitem{Hicken2014}
J.~E. Hicken, Inexact {Hessian}-vector products in reduced-space
  differential-equation constrained optimization, Optimization and Engineering
  15~(3) (2014) 575--608.
\newblock \href {http://dx.doi.org/10.1007/s11081-014-9258-6}
  {\path{doi:10.1007/s11081-014-9258-6}}.

\bibitem{Biros2005a}
G.~Biros, O.~Ghattas, Parallel {Lagrange}-{Newton}-{Krylov}-{Schur} methods for
  {PDE}-constrained optimization. part {I}: The {Krylov}-{Schur} solver, {SIAM}
  Journal on Scientific Computing 27~(2) (2005) 687--713.
\newblock \href {http://dx.doi.org/10.1137/s106482750241565x}
  {\path{doi:10.1137/s106482750241565x}}.

\bibitem{Biros2005}
G.~Biros, O.~Ghattas, Parallel {Lagrange}-{Newton}-{Krylov}-{Schur} methods for
  {PDE}-constrained optimization. part {II}: The {Lagrange}-{Newton} solver and
  its application to optimal control of steady viscous flows, {SIAM} Journal on
  Scientific Computing 27~(2) (2005) 714--739.
\newblock \href {http://dx.doi.org/10.1137/s1064827502415661}
  {\path{doi:10.1137/s1064827502415661}}.

\bibitem{Gill2005}
P.~E. Gill, W.~Murray, M.~A. Saunders, {SNOPT}: an {SQP} algorithm for
  large-scale constrained optimization, SIAM Review 47~(1) (2005) 99--131.
\newblock \href {http://dx.doi.org/10.1137/S0036144504446096}
  {\path{doi:10.1137/S0036144504446096}}.

\bibitem{Schittkowski2009}
K.~Schittkowski,
  \href{https://www.researchgate.net/publication/238690491_NLPQLP_A_Fortran_Implementation_of_a_Sequential_Quadratic_Programming_Algorithm_with_Distributed_and_Non-Monotone_Line_Search_-_User's_Guide_Version_30}{{NLPQLP}:
  a {F}ortran implementation of a sequential quadratic programming algorithm
  with distributed and non-monotone line search - user's guide, version 3.0},
  Report, Department of Mathematics, University of Bayreuth (2009).
\newline\urlprefix\url{https://www.researchgate.net/publication/238690491_NLPQLP_A_Fortran_Implementation_of_a_Sequential_Quadratic_Programming_Algorithm_with_Distributed_and_Non-Monotone_Line_Search_-_User's_Guide_Version_30}

\bibitem{Waechter2005}
A.~Wächter, L.~T. Biegler, On the implementation of an interior-point filter
  line-search algorithm for large-scale nonlinear programming, Mathematical
  Programming 106~(1) (2005) 25--57.
\newblock \href {http://dx.doi.org/10.1007/s10107-004-0559-y}
  {\path{doi:10.1007/s10107-004-0559-y}}.

\bibitem{Hicken2013}
J.~Hicken, J.~Alonso, Comparison of reduced- and full-space algorithms for
  {PDE}-constrained optimization, in: 51st {AIAA} Aerospace Sciences Meeting
  including the New Horizons Forum and Aerospace Exposition, American Institute
  of Aeronautics and Astronautics, 2013.
\newblock \href {http://dx.doi.org/10.2514/6.2013-1043}
  {\path{doi:10.2514/6.2013-1043}}.

\bibitem{rol-website}
D.~Kouri, D.~Ridzal, G.~von Winckel,
  \href{https://trilinos.github.io/rol.html}{{R}apid {O}ptimization {L}ibrary}
  (2020 (acccessed May 22, 2020)).
\newline\urlprefix\url{https://trilinos.github.io/rol.html}

\bibitem{Saad1986}
Y.~Saad, M.~H. Schultz, {GMRES}: a generalized minimal residual algorithm for
  solving nonsymmetric linear systems, SIAM Journal on Scientific and
  Statistical Computing 7~(3) (1986) 856--869.
\newblock \href {http://dx.doi.org/10.1137/0907058}
  {\path{doi:10.1137/0907058}}.

\bibitem{Gauger2008_OptimalDesignOneShot}
N.~R. Gauger, A.~Griewank, J.~Riehme, Extension of fixed point {PDE} solvers
  for optimal design by one-shot method, European Journal of Computational
  Mechanics 17~(1-2) (2008) 87--102.
\newblock \href {http://dx.doi.org/10.3166/remn.17.87-102}
  {\path{doi:10.3166/remn.17.87-102}}.

\bibitem{Kelley1998}
C.~T. Kelley, D.~E. Keyes, Convergence analysis of pseudo-transient
  continuation, {SIAM} Journal on Numerical Analysis 35~(2) (1998) 508--523.
\newblock \href {http://dx.doi.org/10.1137/s0036142996304796}
  {\path{doi:10.1137/s0036142996304796}}.

\bibitem{Broyden1969}
C.~G. Broyden, A new double-rank minimization algorithm, Notices American
  Mathematical Society 16 (1969) 670.

\bibitem{Fletcher1970}
R.~Fletcher, A new approach to variable metric methods, The Computer Journal
  13~(3) (1970) 317--322.

\bibitem{Goldfarb1970}
D.~Goldfarb, A family of variable-metric methods derived by variational means,
  Mathematics of Computation 24~(109) (1970) 23--26.
\newblock \href {http://dx.doi.org/10.2307/2004873}
  {\path{doi:10.2307/2004873}}.

\bibitem{Shanno1970}
D.~F. Shanno, Conditioning of quasi-newton methods for function minimization,
  Mathematics of Computation 24~(111) (1970) 647--657.
\newblock \href {http://dx.doi.org/10.2307/2004840}
  {\path{doi:10.2307/2004840}}.

\bibitem{Nocedal1980}
J.~Nocedal, Updating quasi-newton matrices with limited storage, Mathematics of
  Computation 35~(151) (1980) 773--773.
\newblock \href {http://dx.doi.org/10.1090/s0025-5718-1980-0572855-7}
  {\path{doi:10.1090/s0025-5718-1980-0572855-7}}.

\bibitem{ifpack-website}
M.~Sala, M.~Heroux, \href{https://trilinos.github.io/sacado.html}{{I}fpack}
  (2020 (acccessed May 22, 2020)).
\newline\urlprefix\url{https://trilinos.github.io/sacado.html}

\bibitem{Reed1973}
W.~H. Reed, T.~R. Hill,
  \href{http://www.osti.gov/scitech/servlets/purl/4491151}{Triangular mesh
  methods for the neutron transport equation}, Tech. rep., Los Alamos National
  Laboratory, Los Alamos, New Mexico, USA, lA-UR-73-479 (1973).
\newline\urlprefix\url{http://www.osti.gov/scitech/servlets/purl/4491151}

\bibitem{Roe1981}
P.~L. Roe, Approximate {Riemann} solvers, parameter vectors, and difference
  schemes, Journal of Computational Physics 43~(2) (1981) 357--372.
\newblock \href {http://dx.doi.org/10.1016/0021-9991(81)90128-5}
  {\path{doi:10.1016/0021-9991(81)90128-5}}.

\bibitem{Bangerth2007}
W.~Bangerth, R.~Hartmann, G.~Kanschat, deal.{II}---a general-purpose
  object-oriented finite element library, {ACM} Transactions on Mathematical
  Software 33~(4).
\newblock \href {http://dx.doi.org/10.1145/1268776.1268779}
  {\path{doi:10.1145/1268776.1268779}}.

\bibitem{sacado-website}
D.~M. Gay, E.~T. Phipps,
  \href{https://trilinos.github.io/sacado.html}{{S}acado} (2020 (acccessed May
  22, 2020)).
\newline\urlprefix\url{https://trilinos.github.io/sacado.html}

\bibitem{Jameson1995}
A.~Jameson, Optimum aerodynamic design using {CFD} and control theory, in: 12th
  Computational Fluid Dynamics Conference, American Institute of Aeronautics
  and Astronautics, 1995.
\newblock \href {http://dx.doi.org/10.2514/6.1995-1729}
  {\path{doi:10.2514/6.1995-1729}}.

\bibitem{Sederberg1986}
T.~W. Sederberg, S.~R. Parry, Free-form deformation of solid geometric models,
  SIGGRAPH Computer Graphics 20~(4) (1986) 151--160.
\newblock \href {http://dx.doi.org/10.1145/15886.15903}
  {\path{doi:10.1145/15886.15903}}.

\bibitem{Tezduyar1992}
T.~Tezduyar, M.~Behr, S.~Mittal, A.~Johnson, Computation of Unsteady
  Incompressible Flows with the Stabilized Finite Element Methods: Space-Time
  Formulations, Iterative Strategies and Massively Parallel Implementations,
  Vol. 246, 1992, pp. 7--24.

\bibitem{Truong2008}
A.~H. Truong, C.~A. Oldfield, D.~W. Zingg, Mesh movement for a discrete-adjoint
  newton-krylov algorithm for aerodynamic optimization, {AIAA} Journal 46~(7)
  (2008) 1695--1704.
\newblock \href {http://dx.doi.org/10.2514/1.33836}
  {\path{doi:10.2514/1.33836}}.

\bibitem{Persson2009}
P.-O. Persson, J.~Peraire, Curved mesh generation and mesh refinement using
  {Lagrangian} solid mechanics, in: 47th AIAA Aerospace Sciences Meeting
  including The New Horizons Forum and Aerospace Exposition, American Institute
  of Aeronautics and Astronautics, 2009.
\newblock \href {http://dx.doi.org/10.2514/6.2009-949}
  {\path{doi:10.2514/6.2009-949}}.

\bibitem{Brown2018}
D.~A. Brown, S.~Nadarajah, H.~Yang, P.~Castonguay, H.~Raiesi, K.~Sermeus,
  P.~Germain, Quality-preserving linear elasticity mesh movement algorithm for
  multi-element unstructured meshes, {AIAA} Journal (2018) 1--11\href
  {http://dx.doi.org/10.2514/1.j057463} {\path{doi:10.2514/1.j057463}}.

\bibitem{Orszag1980_SumFact}
S.~A. Orszag, Spectral methods for problems in complex geometries, Journal of
  Computational Physics 37~(1) (1980) 70--92.
\newblock \href {http://dx.doi.org/10.1016/0021-9991(80)90005-4}
  {\path{doi:10.1016/0021-9991(80)90005-4}}.

\bibitem{AndreasGriewank2008}
A.~W. Andreas~Griewank,
  \href{https://www.ebook.de/de/product/23743815/andreas_griewank_andrea_walther_evaluating_derivatives_principles_and_techniques_of_algorithmic_differentiation.html}{Evaluating
  Derivatives: Principles and Techniques of Algorithmic Differentiation},
  CAMBRIDGE UNIV PR, 2008.
\newline\urlprefix\url{https://www.ebook.de/de/product/23743815/andreas_griewank_andrea_walther_evaluating_derivatives_principles_and_techniques_of_algorithmic_differentiation.html}

\bibitem{Pazner2018}
W.~Pazner, P.-O. Persson, Approximate tensor-product preconditioners for very
  high order discontinuous galerkin methods, Journal of Computational Physics
  354 (2018) 344--369.
\newblock \href {http://dx.doi.org/10.1016/j.jcp.2017.10.030}
  {\path{doi:10.1016/j.jcp.2017.10.030}}.

\bibitem{Franciolini2020}
M.~Franciolini, L.~Botti, A.~Colombo, A.~Crivellini, p-multigrid matrix-free
  discontinuous galerkin solution strategies for the under-resolved simulation
  of incompressible turbulent flows, Computers {\&} Fluids 206 (2020) 104558.
\newblock \href {http://dx.doi.org/10.1016/j.compfluid.2020.104558}
  {\path{doi:10.1016/j.compfluid.2020.104558}}.

\bibitem{HOCFD2018}
AIAA, \href{https://how5.cenaero.be/}{5th international workshop on high-order
  {CFD} methods} (2018 (acccessed June 4, 2020)).
\newline\urlprefix\url{https://how5.cenaero.be/}

\bibitem{Carlson2011}
J.-R. Carlson, Inflow/outflow boundary conditions with application to {FUN3D},
  Tech. rep., Langley Research Center, nASA/TM?2011-217181 (2011).

\end{thebibliography}
